\documentclass[a4paper,10pt]{article}

\usepackage{geometry}
\usepackage{lscape}
\usepackage{array}
\usepackage{tikz}
\usepackage{tikz-qtree}

\usepackage{graphicx}
\usepackage{amsmath}
\usepackage{bbm}
\usepackage{cite}
\usepackage[sort&compress,numbers]{natbib}
\usepackage{amsfonts}
\usepackage{amssymb}
\usepackage{amsmath}
\usepackage{latexsym}
\usepackage{color}
\usepackage{url}
\usepackage{ntheorem}
\usepackage{braket}
\usepackage[normalem]{ulem}
\usepackage{enumerate}
\usepackage[colorlinks=true]{hyperref}

\newcommand{\mygb}{\red{b}}

\newcommand{\ominf}{\blue{\omega}}

\newcommand{\myx}{{\red{x}}}
\newcommand{\myxh}{{\red{\vec y}}}

\newcommand{\apar}{\blue{a}}

\newcommand{\sparallel}{{\scalebox{.6}{$\red{\parallel}$}}}
\newcommand{\repsilon}{ {\varepsilon}}

\newcommand{\rhozero}{\red{1}}

\newcommand{\rhozerosquared}{\red{1}}

\newcommand{\rphi}{\red{\psi}}

\newcommand{\nofbox}[1]{#1}
\newcommand{\sfbox}[1]{\fbox{#1}}

\newcommand{\Db}{\mcD}
\newcommand{\Rphi}{\red{\phi}}

\newcommand{\myV}{\mathbb{X}}
\newcommand{\neww}{\red{f_{i,\epsilon}}}

\newcommand{\horo}{\red{\mathfrak{h}}}

\newcommand{\horoG}{\red{\mycal P}}
\newcommand{\Nk}{\red{N^{n-1}}}
\newcommand{\Nkone}{\red{N^{n-1}_1}}
\newcommand{\Nktwo}{\red{N^{n-1}_2}}
\newcommand{\Nki}{\red{N^{n-1}_i}}
\newcommand{\Nka}{\red{N^{n-1}_i}}
\newcommand{\cNb}{\red{\mycal K_b}}
\newcommand{\Ricb}{\red{\mathring{R}}}
\newcommand{\resth}{\red{\mathsf{h}}}
\newcommand{\hk}{\red{h_{k}}}
\newcommand{\hka}{\red{h_{k_i}}}
\newcommand{\hkb}{\red{h_{\bar k}}}

\newcommand{\zMone}{\red{{\hat M}_i}}
\newcommand{\hath}{\red{\hat h }}
\newcommand{\hatR}{\red{\hat R }}

\newcommand{\HMext}{\red{ {M_\ext}}}

\newcommand{\enmom}{\red{m}}
\newcommand{\benmom}{{\mathbf{m}}}

\newcommand{\myone}{\red{1}}
\newcommand{\mytwo}{\red{2}}

\newcommand{\ptcheck}[1]{\ptc{checked on #1}}

\newcommand{\ntheta}{{\psi}}

\newcommand\ben{\begin{enumerate}}
\newcommand\een{\end{enumerate}}
\newcommand\bit{\begin{itemize}}
\newcommand\eit{\end{itemize}}

\newcommand{\qedskip}{\qed\medskip}

\newcommand{\blue}[1]{{\color{blue}#1}}
\newcommand{\red}[1]{{\color{red}#1}}

\newcounter{mnotecount}[section]

\renewcommand{\themnotecount}{\thesection.\arabic{mnotecount}}

\newcommand{\mnote}[1]
{\protect{\stepcounter{mnotecount}}$^{\mbox{\footnotesize
$
\bullet$\themnotecount}}$ \marginpar{
\raggedright\tiny\em
$\!\!\!\!\!\!\,\bullet$\themnotecount: #1} }

\newtheorem{theorem}{\sc  Theorem\rm}[section]

\newtheorem{Theorem}[theorem]{\sc  Theorem\rm}

\newtheorem{conjecture}[theorem]{\sc  Conjecture\rm}

\newtheorem{Corollary}[theorem]{\sc  Corollary\rm}

\newtheorem{Lemma}[theorem]{\sc Lemma\rm}

\newtheorem{Proposition}[theorem]{\sc Proposition\rm}
\newtheorem{remark}[theorem]{\sc Remark\rm}
\newtheorem{Remark}[theorem]{\sc Remark\rm}

\newcommand{\jlcax}[1]{}
\newcommand{\eean}{\nonumber\end{eqnarray}}

\newcommand{\kk}[1]{}

\newcommand{\beq}{\begin{equation}}

\newcommand{\T}{\mathbb T}

\newcommand{\FS}       
                  {F}

\newcommand{\HS} 
       {H_{\mbox{\scriptsize volume}}}

{\ptc{this should be removed in the oberwolfach version}}%

\newcommand{\mc}{\red{m_{\mathrm c}}}

\newcommand{\zD}{\mathring{D}}%

\newcommand{\ourU}{\mathbb U}%
\newcommand{\eeal}[1]{\label{#1}\end{eqnarray}}
\newcommand{\C}{{\mathbb C}}
\newcommand{\bed}{\begin{deqarr}}
\newcommand{\eed}{\end{deqarr}}
\newcommand{\bedl}[1]{\begin{deqarr}\label{#1}}
\newcommand{\eedl}[2]{\arrlabel{#1}\label{#2}\end{deqarr}}

\newcommand{\mcO}{{\mycal O}}
\newcommand{\mcU}{{\mycal U}}

\newcommand{\bel}[1]{\begin{equation}\label{#1}}
\newcommand{\bea}{\begin{eqnarray}}
\newcommand{\bean}{\begin{eqnarray}\nonumber}
\newcommand{\beal}[1]{\begin{eqnarray}\label{#1}}
\newcommand{\eea}{\end{eqnarray}}

\newcommand{\nn}{\nonumber}

\newcommand{\qed}{\hfill $\Box$ \medskip}
\newcommand{\proof}{\noindent {\sc Proof:\ }}
\newcommand{\be}{\begin{equation}}
\newcommand{\eeq}{\end{equation}}
\newcommand{\ee}{\end{equation}}
\newcommand{\beqa}{\begin{eqnarray}}
\newcommand{\eeqa}{\end{eqnarray}}
\newcommand{\beqan}{\begin{eqnarray*}}
\newcommand{\eeqan}{\end{eqnarray*}}
\newcommand{\ba}{\begin{array}}
\newcommand{\ea}{\end{array}}

\newcommand{\mcD}{{\mycal D}}

\newcommand{\warn}[1]
{\protect{\stepcounter{mnotecount}}$^{\mbox{\footnotesize
$
\bullet$\themnotecount}}$ \marginpar{
\raggedright\tiny\em
$\!\!\!\!\!\!\,\bullet$\themnotecount: {\bf Warning:} #1} }

\newcommand{\R}{\mathbb R}

\newcommand{\N}{\mathbb N}
\newcommand{\Z}{\mathbb Z}

\newcommand{\eq}[1]{(\ref{#1})}

\newcommand{\ext}{\mathrm{ext}}

\newcommand{\ptc}[1]{\mnote{{\bf ptc:}#1}}

\newcommand{\beqar}{\begin{deqarr}}
\newcommand{\eeqar}{\end{deqarr}}

\newcommand{\beaa}{\begin{eqnarray*}}
\newcommand{\eeaa}{\end{eqnarray*}}

\newcommand{\hrho}{\hat\rho}

\newcommand{\mv}{\omega}
\newcommand{\mage}[1]{#1}

\DeclareFontFamily{OT1}{rsfs}{}
\DeclareFontShape{OT1}{rsfs}{CGNPm}{n}{ <-7> rsfs5 <7-10> rsfs7 <10-> rsfs10}{}
\DeclareMathAlphabet{\mycal}{OT1}{rsfs}{CGNPm}{n}

{\catcode `\@=11 \global\let\AddToReset=\@addtoreset}
\AddToReset{equation}{section}

{\catcode `\@=11 \global\let\AddToReset=\@addtoreset}
\AddToReset{figure}{section}

{\catcode `\@=11 \global\let\AddToReset=\@addtoreset}
\AddToReset{table}{section}

\renewcommand{\ptcheck}[1]{} 
\renewcommand{\red}[1]{#1}
\renewcommand{\blue}[1]{#1}

\begin{document}
\title{Hyperbolic energy and Maskit gluings\protect\thanks{Preprint UWThPh-2021-25}}

\author{Piotr T. Chru\'{s}ciel\thanks{University of Vienna, Faculty of Physics} \thanks{
{\sc Email} \protect\url{piotr.chrusciel@univie.ac.at}, {\sc URL} \protect\url{homepage.univie.ac.at/piotr.chrusciel}}
\\
{Erwann
Delay}\thanks{ Avignon Universit\'e, Laboratoire de Math\'ematiques d'Avignon,
F-84916 Avignon}
\thanks{Aix Marseille Universit\'e -- F.R.U.M.A.M.- CNRS-- F-13331 Marseille}
 \thanks{{\sc Email} \protect\url{Erwann.Delay@univ-avignon.fr}, {\sc URL} \protect\url{https://erwanndelay.wordpress.com}}
\\
{Raphaela Wutte}\thanks{TU Wien, Faculty of Physics, Institute of Theoretical Physics}
\thanks{
{\sc Email} \protect\url{rwutte@hep.itp.tuwien.ac.at} {} \protect\url{}}
}
\maketitle

\begin{abstract}
We derive a formula for the energy of asymptotically locally hyperbolic
(ALH) manifolds obtained by a gluing at infinity of two ALH manifolds. As
an application we show that there exist three-dimensional conformally compact ALH manifolds
 either without boundary or with toroidal black hole boundary,
 with connected conformal infinity of higher genus, with constant negative scalar curvature, and
with negative mass.
\end{abstract}

\tableofcontents


\section{Introduction}

In~\cite{ILS} Isenberg, Lee and Stavrov, inspired by~\cite{MazzeoPacardMaskit}, have shown how to glue together two
asymptotically locally hyperbolic (ALH) general relativistic initial data
sets by performing a boundary connected sum, which they referred to as a
``Maskit gluing''. The resulting initial data set has a conformal
boundary at infinity which is a connected sum of the original ones. A
variation of this construction has been presented in~\cite{ChDelayExotic}. It
is of interest to analyse the properties of the initial data sets so
obtained. The aim of this work is to address this question in the time
symmetric case, with vanishing extrinsic curvature tensor $K_{ij}$.

We start by a short presentation of  the boundary-gluing construction
of~\cite{ChDelayExotic} in Section~\ref{s2VII21.1}. The construction involves a certain amount of
freedom which we make precise, showing that the gluing results in whole
families of new ALH metrics. As a particular case, in Section~\ref{s29VIII21.1} we apply the construction
to Birmingham-Kottler metrics. This provides  new families of ALH metrics
with apparent-horizon boundaries with more than one component, and with
locally explicit metric when the mass parameter is zero.
One thus obtains time-symmetric vacuum initial data for spacetimes containing
several apparent horizons; such initial data sets \red{are expected to}
evolve to spacetimes with multiple black holes.

Next, an important global invariant of asymptotically hyperbolic general
relativistic initial data sets is the total energy-momentum vector $\benmom \equiv (\enmom_\mu)$
 \cite{ChHerzlich,Wang,ChNagyATMP}
 (compare~\cite{AbbottDeser,ChruscielSimon,deHaro:2000xn})
 when the conformal metric at infinity is that of the round sphere (we talk of \emph{AH metrics} then),
 and the total mass $\enmom$ for the remaining topologies at infinity; this is reviewed in Section~\ref{s3VII21.1}.
It turns out that the formulae for the mass after  the gluings of both~\cite{ILS} and~\cite{ChDelayExotic} are relatively simple. This is analysed in Section~\ref{s5VII21.2}, where we derive
formulae~\eqref{29VII21.1}-\eqref{29VII21.1asdfILS}. This is the first main result of this paper.

A quick glance at \eqref{29VII21.1asdf}
suggests very strongly that the boundary-gluing of (the space-part of) two Horowitz-Myers metrics, which both have negative mass aspect functions, will lead to a manifold with higher-genus boundary at infinity and with \emph{negative mass}.
This turns out, however, to be subtle because of correction terms that are inherent to the constructions.
In fact, negativity of the total mass is far from clear for the Isenberg-Lee-Stavrov gluings, because these authors use the conformal method,
which changes the mass integrand
in a way which appears difficult to control  in the neck region.
Things are clearer when the localised boundary-gluing of~\cite{ChDelayExotic} is used, and in Section~\ref{s29VII21.1} we show that negativity indeed holds.
 We thus construct \emph{three-dimensional conformally compactifiable ALH manifolds with constant scalar curvature,   either without boundaries at finite distance or with a toroidal minimal surface as a boundary at finite distance, with connected boundary at infinity of higher genus topology, and with negative mass.}
  This is the second main result of this paper. As already hinted to, such metrics provide time-symmetric initial data hypersurfaces for vacuum spacetimes with negative cosmological constant. Here one should keep in mind that the Horowitz-Myers metrics have toroidal topology at infinity, and that the  higher-genus Birmingham-Kottler metrics with negative mass are either nakedly singular, or have a totally geodesic boundary with the same genus as the boundary at infinity, or acquire a conformal boundary at infinity with  two components, and contain an apparent horizon, after a doubling across the boundary at finite distance.

It is  conceivable   that a Maskit gluing of a Horowitz-Myers metric with a spherical Kottler (``Schwarzschild-anti de Sitter'') black hole with  small mass will provide an example of a conformally compact vacuum black hole with toroidal infinity and negative mass.   We plan to return to this question in the future.

  To put our work in perspective,  the currently known or conjectured lower bounds for the mass of asymptotically Birmingham-Kottler manifolds are summarised in Table~\ref{T29VIII21.1}. In particular we note the four-dimensional negative-mass ALH metrics of \cite{PedersenEH,Clarkson:2006zk} (compare \cite{ChenZhang,DoldNegativeMass}),
 where conformal infinity is a quotient of a sphere.

 After this work was completed,   large families of ALH metrics with constant negative scalar curvature and with mass of any sign   have been constructed in~\cite{ChDelayNegative} by extending the methods introduced here, with no, or one, or more  black-hole-boundaries with prescribed topology, and with any higher-genus topology at  conformal infinity.

 We end our introductory remarks by noting
  that while there exists a notion of mass for general locally asymptotically hyperbolic metrics with constant scalar curvature which admit asymptotic static potentials~\cite{deHaro:2000xn}, nothing is known about the sign of the mass for those which are not asymptotically Birmingham-Kottler.
  See also \cite{HerzlichMAH} for a discussion of the issues occurring in this context.
  \tikzset{baseline,every tree node/.style={align=center,anchor=north}}

  \begin{table}
    \begin{tikzpicture}[level 1/.style={sibling distance=-35, level distance=40},
      level 2/.style={level distance=40}]
    \Tree[.{Asymptotically Birmingham-Kottler metrics; \red{$m_{\mathrm{crit}}<0$}}
    [.{\fbox{canonical spherical}}
      [.{no bdry}
          {\nofbox{$\ge 0$ \cite{ChDelayHPETv1}}}
       ]
       [.bdry
          {\nofbox{$\ge 0$ \cite{ChGallowayHPET}}}
       ]
    ]
    [.{\qquad \fbox{Ricci flat conf.\ infinity}
       }
    [[ [.{good spin \cite{WangConformal,CMT}}
        [.{no bdry} {$\ge 0$}
        ]
        [.bdry
            [.$\ge 0$ ]
        ]
      ]
      [.{otherwise}
        [.{no bdry}
            {\nofbox{\red{$\ge m_{\mathrm{crit}}$ ? \cite{BCHMN}}}}
        ]
        [.bdry
            {\nofbox{$\, \exists \   m \le 0$ ?}}
        ]
      ]]]
   ]
  [.{\fbox{other conf.\ infinity}}
     [.{no bdry}
        [.{\sfbox{$\exists \ m< 0$ [0]} \cite{PedersenEH}}
          [.{\nofbox{\red{$\ge m_{\mathrm{crit}}$ ??}}}
          ]
        ]
      ]
      [.bdry  [.{$\mu<0$}
          [.{\nofbox{{$\ge m_{\mathrm{crit}}$ \cite{LeeNeves}}}}
          ]
        ]
         [.{otherwise}
          {\nofbox{\red{$\ge m_{\mathrm{crit}}$ ??}}}
        ]
     ]
  ]
  ]
    \end{tikzpicture}
       \caption{Mass inequalities for asymptotically Birmingham-Kottler metrics. A double question mark indicates that no results are available; a single one indicates existence of partial results. The shorthand ``bdry'' refers to a black-hole boundary. ``Good spin'' denotes a topology where the manifold is spin \emph{and} the spin structure admits asymptotic Killing spinors.
       The case ``other conformal infinity'' includes higher genus topologies when the conformal boundary at infinity is two-dimensional, but also e.g.\ quotients of spheres in higher dimensions. Reference [0] is this work.
       Finally, $\mu$ is the mass aspect function. The  critical value of the mass $m_{\mathrm{crit}}$, assuming it exists, is expected to depend upon the conformal structure of the boundary at infinity.
       \label{T29VIII21.1}}
       \end{table}

       \section{Localised boundary-gluing of ALH metrics}
       \label{s2VII21.1}

      Our analysis to follow is motivated by
      the localised boundary-gluings
      of \red{ALH} manifolds, or initial data sets, as in
      \cite[Section~3.5]{ChDelayExotic}. In this section we present a somewhat more general version of
      these gluings.

       Before entering the subject, some comments on our terminology are in order. We  say that a metric $g$ on a manifold without boundary $M$ has a \emph{conformal completion} $(\bar M, \bar g)$  if $\bar M$ is a manifold with boundary such that  $\bar M= M \cup \partial \bar M$, and if
      there exists a function $\Omega\ge 0$ on $\bar M$ which vanishes precisely on $\partial \bar M$, with $d\Omega$ without zeros on $\partial \bar M$, and with $g=\Omega^{-2} \bar g$ on $M$.
      This definition generalises to manifolds $M$ with boundary, in which case $\partial \bar M$ will be the union of the original boundaries of $\partial M$, where $\Omega $ is strictly positive, and the new ones where $\Omega $ vanishes; the new ones are referred to as \emph{boundaries at conformal infinity}. We   say that $(M,g)$ is \emph{conformally compactifiable}  when $\bar M$ is compact. We say that $(M,g)$ is \emph{asymptotically locally hyperbolic}  (ALH) if the scalar curvature of $g$ approaches a constant when the boundary at infinity is approached. We say that an ALH manifold is \emph{asymptotically hyperbolic} (AH) if the conformal class of $\bar g$ on the conformal boundary at infinity is that of a round sphere. \emph{Asymptotically Birmingham-Kottler (ABK) metrics} are defined as  metrics which asymptote to the Birmingham-Kottler metrics of Section~\ref{s29VIII21.1} below. The Birmingham-Kottler metrics themselves are ALH, which can be seen by setting  $\Omega= 1/r$ in \eqref{26XI21.1}, and noting that they have constant scalar curvature since they solve the time-symmetric general relativistic scalar constraint equation with (negative) cosmological constant.

      An ABK metric can equivalently be defined as an ALH metric such that the conformal class of its conformal metric at infinity contains an Einstein metric. Note that this is always the case for three-dimensional manifolds, hence two-dimensional conformal boundaries, by the uniformisation theorem; thus ALH metrics are necessarily ABK in three dimensions, but this is not the case anymore in higher dimensions. For ABK metrics we can introduce Fefferman-Graham coordinates based on the Einstein representative of the conformal metric at infinity. In these coordinates, the asymptotic expansion of $g$ will coincide with that of a BK metric, say $b$, written in the same coordinate system, up to some order.
      This decay order can be measured using $b$-ON frames; equivalently, by measuring the decay of the $b$-norm of $g-b$; the $g$-norm  of $g-b$ would give an equivalent result too. The mass integrals of Section~\ref{ss5VII21.3}  are well defined and convergent if the decay of $g-b$ so understood is,  roughly speaking, $o(r^{-n/2})$.

      The  differentiability requirements of $\bar g$ at the conformal boundary at infinity  often need to be added in the definitions above, and depend upon the problem at hand. Here we will be interested in a class of manifolds with well defined mass, as will be made precise in Section~\ref{ss5VII21.3}.

      The reader is warned that there is no consistency in the literature concerning terminology. While our definition of ALH coincides with that of several authors, some other authors use AH for what we call ALH here. However, we find it natural to reserve the name AH for the special case where the metric is asymptotic to that of hyperbolic space.

 After these preliminaries, we are ready to proceed with our construction.
      We start with points $p_1$, $p_2$, lying on the conformal boundary of two
      \red{ALH} vacuum initial data sets $(M_1,g_1,K_1)$ and $(M_2,g_2,K_2)$  . (An identical construction applies when $p_1$
      and $p_2$ belong to the same manifold; then the construction provides instead
      a handle connecting a neighborhood of $p_1$ with a neighborhood of $p_2$. Instead of vacuum initial data one can also take e.g.\ data satisfying
      the dominant energy condition; the construction will preserve this. Further, if $K_1\equiv 0 \equiv K_2$, then one can have $K\equiv 0$ throughout the construction.)
       We  assume that both $(M_1, g_1,K_1)$ and $(M_2, g_2,K_2)$ have  extrinsic curvature tensors asymptoting to zero, and  have well-defined total energy-momentum, cf.~Section~\ref{s3VII21.1}.
      As shown in \cite{ChDelayExotic} for deformations of data sets preserving the
      vacuum condition
      and for scalar curvature deformations preserving an
      inequality, or in \cite[Appendix~A]{ChDelayHPETv1}
      for deformations of data sets preserving the dominant energy condition,
       for all $\epsilon>0$ sufficiently small we can construct new initial data sets
      $(M_i, g_{i,\red{\epsilon} },K_{i,\red{\epsilon} })$, $i=1,2$, such
      that the metrics coincide with the hyperbolic metric in coordinate half-balls
      $\mcU_ {\myone,\epsilon }$
       of radius $\epsilon$ around $p_1$ and
      $\mcU_ {\mytwo,\epsilon }$ around $p_2$,
       and the $K_{i,\red{\epsilon}
      }$'s are zero there. Here the coordinate half-balls refer to coordinates on
      the upper half-space model of hyperbolic space ${\mathcal H}^n$,  in which we
      have
      \begin{equation}
       \label{2VI21.5}
       {\mathcal H}^n=\R^n_+=\{(\myxh,\myx)\in\R^{n-1}\times (0,\infty)\}
       \,,
      \end{equation}
      with the metric
      \begin{equation}
       \label{2VI21.6}
      b=\frac{|d \myxh|^2+d\myx^2}{\myx^2}
      \,.
      \end{equation}

      In order to avoid a proliferation of indices we now choose  $\epsilon$
        so that the deformation described above has been carried out, with
      corresponding ALH metrics $(M_i, g_{i,\red{\epsilon}
      },K_{i,\red{\epsilon} })$, $i=1,2$, and  coordinate half-balls
      $\mcU_ {i,\epsilon }$, which we write from now on as $(M_i, g_{i},K_{i})$,
       and   $\mcU_{i}$. It should be kept in mind
      that different values of $\epsilon$ lead to different initial data sets near the
      gluing region, but the data sets remain the original ones, hence identical, away from the gluing region,  which can be chosen as small as desired.

      The above shows in which sense the exterior curvature tensor $K$ is
      irrelevant for the problem at hand. Therefore, from now on, for simplicity we will only
      consider initial data sets with $K\equiv 0$ in the gluing region.

\begin{figure}
        \centering
       \phantom{xx}\includegraphics[width=.9\textwidth]{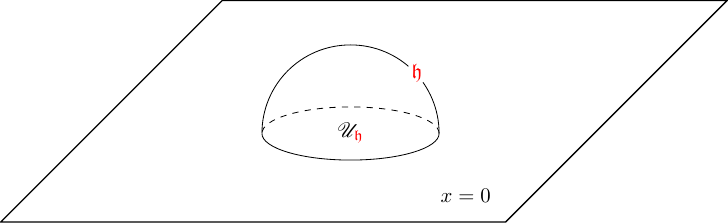}
          \caption{The ``thin component'' $\mcU_\horo$ and its boundary $\horo$ in the half-space model, where hyperbolic space is represented as a half-space $\{\myx >0\}$, and the conformal boundary at infinity is the hyperplane $\{\myx =0\}$. After an  ``exotic gluing'' has been performed, the metric becomes exactly hyperbolic inside $\mcU_\horo$.}
        \label{F2VII21.1}
      \end{figure}

      The gluing construction uses  \emph{hyperbolic hyperplanes} $\horo$ in the
      hyperbolic-space region of $(M_i,g_i)$. These are defined, in the half-space model, as
      half-spheres with centres on the hyperplane $\{\myx=0\}$ in the coordinates of
      \eqref{2VI21.6} and which, for our purposes, are entirely contained in the half-balls $\mcU_i$. The conditionally
      compact, in $\R^n$, component of ${\mathcal H}^n$ separated by $\horo$ will
      be referred to as the \emph{thin component}, denoted by $\mcU_\horo$, and the
      remaining one will be called the \emph{fat component}; see
      Figure~\ref{F2VII21.1}.

      In what follows we will also invoke the \emph{Poincar\'e
      ball} model, which represents the  $n$-dimensional hyperbolic space
      ${\mathcal H}^n$ as the open unit ball $B^n$ endowed with the metric
      \begin{equation}\label{3XI18.4}
        b = \frac{4}{(1-\red{|\vec x|}^2)^2} \delta
        \,,
      \end{equation}
      where $\delta$ is the Euclidean metric.

      A basic fact is that for every  hyperbolic hyperplane $\horo$ as above
      there exist two isometries of the hyperbolic space, which we denote by $\Lambda_{\horo,\pm}$, such
      that $\Lambda_{\horo,+}$ maps $\horo$  to the
      equatorial hyperplane  of the
      Poincar\'e ball, with the fat region being mapped to the upper hemisphere, while $\Lambda_{\horo,-}$ again maps the  hyperbolic hyperplane to the
      equatorial hyperplane but maps the fat region into the lower
      hemisphere. Using physics terminology, an example  of $\Lambda_{\horo,+}$ is
      provided by a boost along the axis passing through the origin of
      the Poincar\'e ball and the barycenter of the hyperbolic hyperplane.  Given
      $\Lambda_{\horo,+}$, a map $\Lambda_{\horo,-}$ can be obtained by applying to
      $\Lambda_{\horo,+}$ a rotation by $\pi$ around any axis lying on the
      equatorial plane.

      It should be clear that there are many such pairs $ \Lambda_{\horo,\pm}$:
      consider, e.g., the isometries $R_\pm \Lambda_{\horo,\pm}$, where the
      $R_\pm$'s are rotations along the axis joining the north pole and the south pole.

      Let $\horoG$ denote the collection of pairs of isometries
      $(\Lambda_1,\Lambda_2)$ of hyperbolic space with the following property:
      There exist  hyperbolic hyperplanes $\horo_i\subset \mathring\mcU_i$, $i=1,2$,
      such that $\Lambda_1$ is a $\Lambda_{\blue{\horo_{1}},+}$
      and $\Lambda_2$ is a
      $\Lambda_{\horo_2,-}$. Here $ \mathring\mcU_i$ denotes the interior of
      $\mcU_i$.

      For each such pair of  hyperbolic hyperplanes $\horo_i$ the manifolds
      $M_i\setminus \mcU_{\horo_i}$ are manifolds with a non-compact boundary component
      $ \partial \mcU_{\horo_i} = \horo_i$, extending to the conformal boundary at infinity
      of $M_i$, with the hyperbolic metric near the boundary $\horo_i$.

      Given $(\Lambda_1,\Lambda_2)\in \horoG$ we construct a boundary-glued
      manifold $M_{\Lambda_1,\Lambda_2}$ by gluing the boundaries $\horo_i$ as
      follows: We map the thin complement $\mcU_{\blue{\horo_{1}}}$ of $\blue{\horo_{1}}$
      to the
      lower half of the Poincar\'e ball using $\Lambda_1$. We map the thin
      complement $\mcU_{\horo_2}$ of $\horo_2$ to the upper half of the Poincar\'e
      ball using $\Lambda_2$. We then identify the two manifolds with boundary
      $M_i\setminus \mcU_{\horo_i}$
      along the equatorial plane of the Poincar\'e
      ball using the identity map.

      The metrics on $M_i\setminus \mcU_{\horo_i}$ coincide with the original ones
      (one can think of the maps $\Lambda_{\horo_i}$ as changes of coordinates),
      hence are ALH there by hypothesis. Both metrics are exactly hyperbolic at both
      sides of  the gluing boundary, namely the equatorial plane of the Poincar\'e
      ball, and extend smoothly there. Hence for every pair
      $(\Lambda_1,\Lambda_2)\in \horoG$ the manifold $M_{\Lambda_1,\Lambda_2}$ is a
      smooth ALH manifold.

      \section{Boundary-gluing of Birmingham-Kottler solutions}
      \label{s29VIII21.1}

      An obvious candidate to which our construction can be applied is the space-part of the Birmingham-Kottler (BK) metrics,
      \begin{equation}\label{26XI21.1}
        g = f^{-1} dr^2 + r^2 \hk
        \,,
        \qquad
        f = r^2 + k - \frac{2\mc }{r^{n-2}}
         \,,
      \end{equation}
      where $\mc$ is a constant which will be referred to as the \emph{coordinate mass parameter}, and
      where $\hk $ is an Einstein metric on an $(n-1)$-dimensional manifold $\Nk$
      with scalar curvature equal to
      $$
       R(\hk)= k (n-1)(n-2)
       \,,
       \qquad k \in \{0,\pm1\}
       \,.
      $$
      The associated spacetime
      metric
      $$
       - f dt^2 + g
      $$
      is a solution of the (Lorentzian) vacuum Einstein equations with a negative
      cosmological constant.

      The collection of $(\Nk,\hk)$'s is quite rich:  one should, e.g., keep in mind the existence of many Einstein metrics on higher dimensional spheres, including exotic ones~\cite{BoyerGalickiKollar}.

      We will assume that $\mc$ is in a range so that $f$ has positive zeros, with the largest
      one, denoted by $r_0$, of first order.  The metric is then smooth on the
      product manifold $[r_0,\infty)\times \Nk$, with a totally geodesic boundary
      at $r=r_0$.

      The construction of  Section~\ref{s2VII21.1} applied to  two such manifolds,
      $$
       (M_i= [r_i,\infty)\times \Nki,g_i)
        \,,
      $$
      leads to manifolds with a conformal boundary at infinity with connected-sum topology $\Nkone
      \#\Nktwo $ and a totally geodesic boundary which has two connected
      components, one diffeomorphic to $\Nkone $ and the second to $\Nktwo $.

      One can double each of the original manifolds across their totally geodesic
      boundaries, in which case the doubled manifolds have no boundaries but each
      has a conformal infinity with two components. Performing our construction on
      a chosen pair of boundaries at infinity one obtains an ALH manifold with
      three boundary components, one with topology $\Nkone $, one with topology
      $\Nktwo $, and a third one with topology  $\Nkone \#\Nktwo $.

      One can iterate the construction, obtaining ALH manifolds with an arbitrary
      number of components of the boundary at infinity, and an  arbitrary number
      of totally geodesic compact
      boundary components. The maximal globally hyperbolic development of the resulting time-symmetric general relativistic vacuum Cauchy data will have a Killing vector field defined in a neighborhood of each such boundary (but not globally in general), which becomes the bifurcation surface of a bifurcate Killing horizon for this vector field.
      The case $\mc=0$ is of special interest here, as then the metric is everywhere exactly  hyperbolic locally,  so
      that the construction of~\cite{ChDelayAH}, which glues-in a hyperbolic
      half-ball, is trivial.
      All the boundary gluing constructions described here apply without further due to this case.
      It is likely that some of the metrics constructed above evolve to the spacetimes considered in~\cite{Peldan1,Peldan2}, it might be of some interest to  explore this.

      \section{Hyperbolic mass}
      \label{s3VII21.1}

     \subsection{The definition}
      \label{ss5VII21.3}

     We recall the definition of hyperbolic mass from \cite{ChHerzlich}. We
     consider a family of Riemannian metrics $g$ which approach a ``background
         metric'' $b$ with constant scalar curvature  $-n(n-1)$ as the conformal boundary at infinity is approached.
         We assume
     that $b$ is equipped with a non-empty set of solutions of the \emph{static Killing Initial Data (KID)}
     equations:
     %
     \begin{eqnarray}
      &  \label{eq:1}
       \red{\mathring \Delta}  V -nV =0\,,
      &
     \\
      &
       \label{eq:2} \zD_i\zD_j V = V( \Ricb_{ij} + n
     b_{ij})
      \,,
      &
     \end{eqnarray}
     where   $\Ricb_{ij}$ denotes the Ricci tensor of the
     metric $b$, $\zD$ the Levi-Civita connection of $b$, the operator
     $\red{\mathring \Delta} :=b^{k\ell}\zD_k\zD_\ell $ is the Laplacian of $b$, and we use $D$ for the Levi-Civita connection of $g$. Nontrivial
     triples $(M,b,V)$, where $V$ solves \eq{eq:1}-\eq{eq:2}, are called
     \emph{static Killing Initial Data}.
     Strictly speaking, for the purpose of defining the mass, only the geometry of a neighborhood of the conformal boundary at infinity of  $M$ is relevant. A couple $(b,V)$ will also be called a static KID when $M$ or its  conformal boundary are implicitly understood.
     The functions $V$ will be
     interchangeably referred to as \emph{static potentials} or  static KIDs.
     The (finite dimensional) vector space of static
     potentials of $(M,b)$ will be denoted by $\cNb$.
     All static potentials for the BK  and Horowitz-Myers metrics  are derived in Appendix~\ref{A17VIII21.1}.

     Under conditions on the convergence of $g$ to $b$ spelled-out in \eqref{m1}-\eqref{m5} below, well defined
     global geometric invariants can be extracted from the integrals~\cite{ChHerzlich}
      \be \label{mi12I}
     H(V,b):=\lim_{R\to\infty} \int_{r=R} \ourU^i(V) \red{dS}_i \ee
     where $V\in
     \cNb$, with
     \begin{eqnarray} \label{eq:3.312I} & {}\ourU^i (V):=  2\sqrt{\det
     g}\;\left(Vg^{i[k} g^{j]l} \zD_j g_{kl}
     +D^{[i}V 
     g^{j]k} (g_{jk}-b_{jk})\right)
     \,,
     \end{eqnarray}
     and
     where $dS_i$ are the hypersurface  forms $\partial_i\rfloor dx^1\wedge \cdots \wedge dx^n$.

     From now on, until explicitly indicated otherwise the background metric  $b$ will be a BK metric with $\mc=0$.
      In the coordinates of
     \eqref{26XI21.1} we set
     \begin{equation}\label{3VII21.2}
     \HMext := [R,\infty)\times \Nk
      \,,
     \end{equation}
     for some large $R\in \R^+$, and we consider the following orthonormal frame
     $\{f_i\}_{i=1,n}$ on $\HMext$:
     \be\label{m1}  f_i = r^{-1}\epsilon_i\,, \quad
     i=1,\ldots,n-1\,,\quad
      f_n = \sqrt{r^2+k} \;\partial_r\,, \ee
     where the $\epsilon_i$'s form an orthonormal frame for the metric $\hk $. We
     set
     \be \label{m2} g_{ij}:=g(f_i,f_j)
      \,,
      \quad b_{ij}:=b(f_i,f_j)
      = \left\{
          \begin{array}{ll}
            1, & \hbox{$i=j$;} \\
            0, & \hbox{otherwise.}
          \end{array}
        \right.
     \ee

 In order to avoid a discussion of the special case $n=2$, which is irrelevant for our main purposes here, from now on  in the main body of the paper  we assume that the space-dimension $n$ is larger than or equal to three.

     Assuming that
     \begin{subequations}
      \label{Hm3}
     \begin{gather}\label{Hm3a}
      \displaystyle
      \int_{\HMext} \big( \sum_{i,j} |g_{ij}-b_{ij}|^2 + \sum_{i,
      j,k} |f_k(g_{ij})|^2 \big)r
      \;d\mu_g<\infty\,,
     \\
      \displaystyle
      \int_{\HMext} |R_g-R_b|\;r\;d\mu_g<\infty\,,\label{Hm3b}
    \end{gather}
  \end{subequations}
     \begin{equation}
     \label{m0} \exists \ C > 0 \ \textrm{ such that }\
     C^{-1}b(X,X)\le g(X,X)\le Cb(X,X)\,,
     \end{equation}
     where
     $d \mu_{g}$ denotes the measure associated with
     the metric $g$,
     one finds that the limit in \eqref{mi12I} exists and is finite. If moreover
     one requires that
      \be
       \label{m5}
      \sum_{i,j}
         |g_{ij}-b_{ij}| + \sum_{i, j,k} |f_k(g_{ij})|=
         o(r^{-n/2})\,,
     \ee
     then the ``mass integrals'' \eqref{mi12I} are well defined, in the following sense: Consider
     any two background metrics $b_i$, $i=1,2$, of the form \eqref{26XI21.1},
     \emph{with the same boundary manifold $(\Nk,\hk)$ with $\mc=0$,} in particular with the same value of $k$. Assume
      that $g$ approaches both $b_1$ and $b_2$ at the rates presented above. Then
      there exists an isometry $\Phi$ of $b_1$  such that
     \begin{equation}\label{Htransf} H (V,b_1)=
     H(V\circ \Phi ,b_2)
     \,.
     \ee

     Consider a background $b$ of the form \eqref{26XI21.1}.
     The function
     \begin{equation}\label{20III18.1}
      V_{(0)}(r)=  \sqrt{r^2+k}
     \end{equation}
     satisfies the static KID equations  \eqref{eq:1}.
     Assuming \eqref{Hm3}-\eq{m5}, a somewhat lengthy calculation
     shows that the mass
     integral
     $$
      H(V_{(0)}):=H(\sqrt{r^2+k },b)
     $$
     equals~\cite{ChruscielSimon}
      \begin{eqnarray}
      \nn
      \lefteqn{  \displaystyle H (V_{(0)})
      =
      \lim_{R\to\infty} R^{n-1} (R^2+k )\times
     }
      &&
     \\   &&
        \displaystyle \int_{\{r=R\}} \left(-\sum_{i=1}^{n-1}\left\{\frac {
               \partial e_{ii}}{\partial r}+ \frac {
               k  e_{ii}}{ r(r^2+k )}\right\}+\frac {(n-1)e_{nn}}{r}
        \right) d^{n-1}\mu_{\hk}\,.
     \phantom{xxx}
        \label{massequation1}
     \end{eqnarray}
     Here
     $$
      e_{ij}:=g_{ij}-b_{ij}
     $$
     denotes the frame components  in a $b$-ON frame, with the $n$'th component corresponding to the direction orthogonal to the level sets of $r$.

     In spacetime dimension $3+1$, under the decay conditions spelled above  this simplifies to
      \begin{eqnarray}
       \displaystyle H (V_{(0)})
      =
      \lim_{R\to\infty} R^4
        \displaystyle \int_{\{r=R\}} \left( \frac {2e_{33}}{r}
      - \sum_{i=1}^{2} \frac {
               \partial e_{ii}}{\partial r}
        \right) d^{2}\mu_{\hk}\,.
     \phantom{xxx}
        \label{massequation13d}
     \end{eqnarray}

     A class of ALH metrics of interest are these for which
     \begin{equation}\label{29VII21.3}
       e_{ij} = r^{-n} \mu_{ij} + o (r^{-n})
     \,,
     \end{equation}
     where the $\mu_{ij}$'s depend only upon the coordinates $x^A$ on $\partial M$. One can further specialise the coordinates so that $\mu_{nn}\equiv 0$, but this choice might be unnecessarily restrictive for some calculations. Metrics with this asymptotics are dense in the space of all AH metrics in suitable circumstances~\cite[Theorem~5.3]{DahlSakovich}. The tensor $\mu_{ij}$ will be referred to as the \emph{mass aspect tensor}.

     For metrics in which \eqref{29VII21.3} holds, Equation~\eqref{massequation1} reads
           %
      \begin{eqnarray}
      \displaystyle H (V_{(0)})
      =  \int_{ \partial M}
     \Big(
       (n-1)\mu_{nn}
     +  n \sum_{i=1}^{n-1}\mu_{ii}
        \Big) d^{n-1}\mu_{\hk}
        =:  \int_{ \partial M}
        \mu \,d^{n-1}\mu_{\hk}\,,
        \label{29VII21.4}
     \end{eqnarray}
     where $\mu$ is the \emph{mass aspect function} mentioned in Table \ref{T29VIII21.1}.

     An elegant formula for mass can be derived by  integration by parts in \eqref{mi12I}-{\eqref{eq:3.312I}, leading to the following:
     To every KID $(b,V)$ and conformal-boundary component $\partial M$ one associates a mass
     $$m=m(V)=m(V,\partial M)\equiv   H(V,b)
     $$
     by the formula
       \cite{HerzlichRicciMass} (compare \cite[Equation~(IV.40)]{BCHKK})
     \begin{eqnarray}
      m(V,\partial M)
      &
      =
     &
        - 
        \lim_{x\rightarrow0}\int_{\{x\} \times \partial M}   D^j V
         ( R{}^i{}_j - \frac {R{}}{n}\delta^i_j)
         \,
         d\sigma_i
          \,,
                \label{18IX20.4}
         \end{eqnarray}
     where
     $R_{ij}$ is the Ricci tensor of the metric $g$, $R$ its trace, and
     we have ignored an overall dimension-dependent positive multiplicative factor which is typically included in the physics literature. Here $\partial M $ is a component of conformal infinity,  $x$ is a coordinate near $\partial M$ so that $\partial M$ is given by the equation $\{x=0\}$, and $d\sigma_i := \sqrt{\det g} dS_i$.

     As a special case, consider a triple $(M,g,\red{\hat V})$   satisfying \eqref{eq:1}-\eqref{eq:2}; thus both $(M,g,\red{\hat V})$ and $(M,b,V)$ are static KIDs. Assume that $g$ asymptotes to $b$ as above, and that $\red{\hat V}$ asymptotes to $V$.
     One checks that $V$ can be replaced by $\red{\hat V}$ in  the integrand of \eqref{18IX20.4}, so that
     \begin{eqnarray}
       D^j\red{\red{\hat V}} (R^i{}_j - \frac{R}{n} \delta^i_j)
        &=&
        \frac{1}{\red{\red{\hat V}}} D^j \red{\red{\hat V}}
         \big(
          D^i D_j \red{\red{\hat V}} - \frac{1}{n}\Delta \red{\red{\hat V}} \delta^i_j
          \big)
          \nonumber
     \\
       & = &
        \frac{1}{\red{\red{\hat V}}} D^j \red{\red{\hat V}}
        \big(
         D^i D_j \red{\red{\hat V}} -  \red{\red{\hat V}}  \delta^i_j
         \big)
          \nonumber
     \\
       & = &
        \frac{1}{2\red{\red{\hat V}}} D^i
         (|d
          \red{\red{\hat V}}|^2)
       -   D^i \red{\red{\hat V}}
        \,.
         \label{22VII21.1}
     \end{eqnarray}
     Letting $r=1/x$, \eqref{18IX20.4} becomes
     \begin{eqnarray}
      m(\red{\hat V})
         &= &
        \lim_{R\rightarrow \infty}\int_{r=R}
        \big(  D^i \red{\hat V}
        -
        \frac{1}{2\red{\hat V}} D^i
          (|d
          \red{\hat V}|^2)
       \big)
         \,
         d\sigma_i
          \,.
                \label{22VII21.1a}
         \end{eqnarray}

     As an application, and for further use, we apply \eqref{22VII21.1a} to a Horowitz-Myers metric,
     \begin{equation}\label{3VIII21.1}
       ds^2 = -r^2 dt^2 +
        \underbrace{
          \frac{dr^2}{r^2 - \frac{2\mc }{r^{n-2}}}
      + (r^2 - \frac{2\mc }{r^{n-2}}) d\ntheta^2 + r^2 \resth_{ab}d\theta^a d\theta^b
      }_{=:g}
      \,,
     \end{equation}
     where  $\resth_{ab}d\theta^ad\theta^b$, $a,b\in\{1,\ldots n-2\}$ is a flat $(n-2)$-dimensional metric, $\mc >0$ is a constant and $\ntheta$ has a suitable period to guarantee smoothness at $ r^n=2\mc  $. The background is taken to be $g$ with $\mc =0$, thus $\red{\hat V}= r$, $|d\red{\hat V}|^2 = (r^2 - \frac{2\mc }{r^{n-2}})$, which gives an integrand in \eqref{22VII21.1a} equal to
     \begin{eqnarray}
       \big(  D^r \red{\hat V}
        -
        \frac{1}{2\red{\hat V}} D^r
          (|d
          \red{\hat V}|^2)
       \big)
         \,
          r^{ {n-2}}
          \sqrt{\det \resth}
       =  -  (n-2) \mc  \big(1- \frac{2\mc }{r^{n}}
       \big)
          \sqrt{\det \resth}
          \,,
     \label{3VIII21.2}
     \end{eqnarray}
     and thus total mass
     \begin{equation}\label{3VIII21.11}
       m = -  (n-2)\mc
       \,,
     \end{equation}
     where we have assumed that the area of the conformal boundary at infinity in the metric $d\ntheta^2 + \resth_{ab}dx^a dx^b$ has been normalised to $1$.

     Applying~\eqref{22VII21.1a} to the toroidal Birmingham-Kottler metric
     \begin{equation}\label{3VIII21.12}
       ds^2 = -
       \left(
        r^2 - \frac{2\mc }{r^{n-2}}
        \right)
        dt^2 +
         \underbrace{
          \frac{dr^2}{r^2 - \frac{2\mc }{r^{n-2}}}
      +
       r^2 \resth_{AB}d\theta^A d\theta^B
       }_{=:g}
      \,,
     \end{equation}
     which has $\hat V = \sqrt{r^2 - \frac{2\mc }{r^{n-2}}}$, gives
     \begin{eqnarray}
       \big(  D^r \red{\hat V}
        -
        \frac{1}{2\red{\hat V}} D^r
          (|d
          \red{\hat V}|^2)
       \big)
         \,
           \frac{r^{n-1}}{\hat V}
       =
       \mc (n-2) (n-1)\left(1 + \frac{\mc (n-2)} {r^{n} }\right)
          \,,
     \label{25V22.98}
     \end{eqnarray}
     and
     \begin{equation}\label{25V22.99}
       m = \left(n-1\right)\left(n-2\right)\mc
        \,.
     \end{equation}

     \subsection{The spherical case}
      \label{ss19VII21.1}

     The question arises, what happens with the mass under the boundary gluings of Section~\ref{s2VII21.1}. The case when both manifolds have conformal boundary with spherical topology and with a metric conformal to the standard round metric is simplest, because then the maps
     $(\Lambda_1,\Lambda_2)\in \horoG$ act globally on collar neighbourhoods of $\Nkone$ and $\Nktwo$ in $M_1$ and $M_2$. The initial energy-momenta $\benmom_1$ of $\Nkone$ and $\benmom_2$ of $\Nktwo$ are transformed to $\Lambda_1 \benmom_1$ of
     $\Lambda_2 \benmom_2$. Since the mass integrands are zero in the neck region, where the metric is exactly hyperbolic, one finds that the energy-momentum  $\benmom$ of the boundary-glued manifold is additive:
     \begin{equation}\label{19VII21.1}
       \benmom = \Lambda_1 \benmom_1 + \Lambda_2 \benmom_2
      \,.
     \end{equation}
     This is true for all dimensions $n \ge 3$.

     A more detailed presentation of the spherical case can be found in~\cite{ChDelayHPETv1}, and we note that most of the work there arose from the necessity to control the direction of the momenta $\benmom_a$ under the gluing-in of the exactly-hyperbolic region.

We take this opportunity to comment on the approach to the positive energy theorem in \cite{ChDelayHPETv1}.  The argument did not cover the case $n=3$ because the estimates on the perturbation of the energy-momentum of the metric deformed as in \cite{ChDelayExotic} were not sufficiently strong to conclude. In dimension $n=4$ the argument there applied at the price of a very careful examination of the error terms.
Proposition~\ref{P21V22.1} below, with $\psi_i$ independent of $i$, leads to the following improvement of the analysis in  \cite[Remark~2.1]{ChDelayHPETv1}:

\begin{Corollary}
  \label{C21V22.2}
  Assume that the metric $g$ has the asymptotic behaviour \eqref{21V22.3} with $\sigma>n-1/2$.
In dimensions $n=3,4$ the error terms in~\cite[Equation~(2.4)]{ChDelayHPETv1} can be improved to
\begin{equation}\label{21V22.41}
  m_\mu^{i,\epsilon} -  m_\mu^{i}  = O(\epsilon^{\sigma-1})
  \,,
\end{equation}
in particular the no-spacelike-energy-momentum arguments there also apply in dimension $n=3$, and without the need of the supplementary careful analysis of the error terms in dimension $n=4$ as done in~\cite{ChDelayHPETv1}.
\end{Corollary}

\subsection{Conformal rescaling of \texorpdfstring{$\hk$}{hk}}

The remaining cases,
as well as the Isenberg-Lee-Stavrov gluings,
 require more effort. Some preliminary remarks are in order.

After the gluing has been done,  the initial ``boundary metric'' $\hk $ of \eqref{26XI21.1} will
have to be conformally rescaled in general. Thus, there will be  a function $\rphi>0$, defined on a subset of the boundary at infinity, such that
\begin{equation}\label{5VII21.1}
  \hk = \rphi^2 \hkb
 \,,
\end{equation}
for some constant $\bar k \in \{0,\pm 1\}$.

The following can be justified by considerations revolving around the Yamabe problem and the Obata equation.

First, if $(\Nk,\hk)$ is a round sphere, there exist globally defined non-trivial such conformal factors with $k=\bar k=1$, and they all arise from conformal isometries of the sphere.

Next, when $(\Nk,\hk)$ is a closed manifold with $k=0$ and \eqref{5VII21.1} holds globally, then $\bar k $ must also be zero and $\rphi$ must be constant. Of course any constant will do, and there does not seem to be a geometrically preferred value for this constant.
When calculating the mass of the Horowitz-Myers metrics
 we will normalise the volume of $(\Nk,\hk)$ to one, which removes the ambiguity.

In all other cases $\rphi\equiv 1$ and $k=\bar k$ is the only possibility for globally defined functions $\rphi$  on closed manifolds.

However, we will also need \eqref{5VII21.1} on subsets of $\Nk$. Then non-trivial functions $\rphi$ are possible with $k\ne \bar k$ in dimension $n=3$. By way of example, let
\begin{equation}\label{5VII21.2}
N^{2} = \T^2 \# \T^2
 \,,
\end{equation}
with the gluing occurring across a closed geodesic as in  Figure~\ref{F5V21.1}.
 \begin{figure}
	\centering
 \phantom{xx}\includegraphics{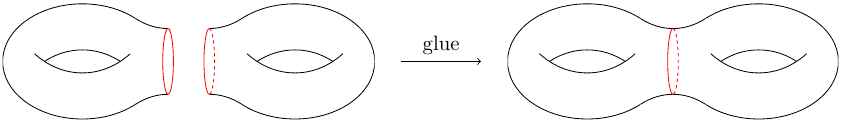}
    \caption{Gluing two solid tori along a disc, the boundary of which becomes a closed geodesic on the boundary of the glued manifold.}
	\label{F5V21.1}
\end{figure}
If we cut the connecting neck along the closed geodesic, each of the factors in \eqref{5VII21.2} can be viewed as a torus with an open disc removed,  $\T^2 \setminus D^2$, and there we can write
  ($k = 0$, $\bar{k} = -1$)
%
\begin{equation}\label{5VII21.2qw}
  h_{-1} = e^{  \mv} h_{0}
 \,,
\end{equation}
where $\mv$ is a solution of the two-dimensional Yamabe equation,
\begin{equation}\label{5VII21.3}
  \Delta_{h_0} \mv =  2 e^{ \mv}
 \,,
\end{equation}
with vanishing Neumann data at $\partial D^2$.

Now, our definition of mass requires that the glued manifold admits non-trivial static potentials for the asymptotic background. This is guaranteed when both the initial manifolds $(\Nka,\hka)$, $i=1,2$, and the   new metric on $\Nkone \# \Nktwo$   are Einstein manifolds. The gluing construction itself does not care about this, but the existence of a static potential on a collar neighborhood of the glued manifold is not clear,
except when the boundary is  two dimensional,
or when at least one   boundary metric is a round sphere  in all dimensions.
The existence of further higher dimensional such examples is unlikely, see~\cite{Bryant}.
%
%
%

Here one should keep in mind that there exist large families of \emph{smoothly compactifiable stationary
vacuum} solutions $(M,g)$
of (Riemannian) Einstein equations such that $(\Nk,\hk)$ is \emph{not} Einstein~\cite{ChDelayKlinger,ACD2,ACD,ChDelayKlingerBH,ChDelayStationary}. Each such metric comes equipped with a definition of mass~\cite{deHaro:2000xn} for metrics with the same conformal infinity. Whether or not our weighted addition formula \eqref{29VII21.1} below applies in these more general circumstances is not clear; we plan to address this in the future.

Returning to our main line of thought, suppose thus that \eqref{5VII21.1} holds on a subset of $\Nk$.
We extend $\rphi$ to a neighbourhood of the conformal boundary by requiring $\partial_r \rphi \equiv 0$.
Let us denote by $(\bar r,\bar x^A)$ coordinates such that the new background metric $\hkb$ takes the form
\begin{equation}\label{26XI21.1zd}
  \bar b:  = \frac{d\bar r^2} {\bar r^2 + \bar k }+ \bar r^2 \hkb
   \,.
\end{equation}
%
We can write
%
\begin{eqnarray}\label{26XI21.1as}
    b   = \frac{d  r^2} {  r^2 +   k }+   r^2 \hk
    = \frac{d  r^2} {  r^2 +   k }+   r^2 \rphi^2 \hkb =
\frac{d\bar r^2} {\bar r^2 + \bar k }+ \bar r^2
\red{
( \hkb
  + \delta h )
}
  \,,
\end{eqnarray}
 {with}
$$
 \delta h(\partial_r, \cdot) = 0
  \,.
$$
Let us denote by $x^A$ the local coordinates at $\partial M$.
As shown in Appendix~\ref{sApp23IX21.1} below (cf.\ \eqref{expxAappendix} and \eqref{exprappendix})
 we find for large $r$
%
\begin{align}
 \bar{x}^A ={}& x^A +\frac{\Db^A \rphi }{2 \rphi r^2}  + O(r^{-4})
      \,,  \label{20IX21.11}\\
\bar{r} ={}& \rphi r \left(1 +
\frac{\left(-\Db^A \Db_A\red{\log}  \rphi + (n-2){\Db}^A \red{\log}\rphi  {\Db}_A \red{\log}\rphi\right) }{2(n-1) r^2}
+ O(r^{-4}) \right) \label{20IX21.11a}
 \,,
\end{align}
where $\Db$ is the Levi-Civita connection of $\hk$.
The pattern is, that the expansion jumps by two powers of $r$ until the threshold associated with the mass aspect tensor.

It follows from \eqref{20IX21.11}-\eqref{20IX21.11a} (see \eqref{conditionaBC}, Appendix~\ref{A17VIII21.2}) that
an asymptotic expansion of the initial metric of the form
\begin{equation}\label{5VII21.5}
  g(f_i,f_j)- b_{ij}=  O(r^{-n})
\end{equation}
(which is observed in BK metrics)
will \emph{not} be preserved by the above transformations, except possibly 1) in dimension $n=3$ or 2) when $k= \bar k $ and the conformal factor arises from a conformal isometry of $h_k$.
This is addressed in Appendix~\ref{A17VIII21.2}. In particular it is shown there, in space-dimension three and in the gauge $\mu_{ri}=0$,  that the mass aspect tensor $\mu_{ij}$ of \eqref{29VII21.3} transforms as
\begin{equation}\label{20IX21.22}
  \mu_{AB} \mapsto \bar \mu_{AB} = \rphi \mu_{AB}
  \,.
\end{equation}

More generally, one can define the mass-aspect tensor as the coefficient of $r^{-n}$ in the asymptotic expansion of the frame components of the metric. Then \eqref{20IX21.22} remains valid in odd space dimensions because of the structure \eqref{20IX21.11}-\eqref{20IX21.11a} of the expansion of the coordinate transformation, with powers jumping by two. But \eqref{20IX21.22} does not hold in even space dimensions: in Appendix~\ref{ss20IX21.1} we derive the transformation formula in dimension $n=4$.

\subsection{Gluing in three space dimensions}
 \label{s5VII21.2}

We consider
 the mass of the manifold obtained by boundary-gluing  two three-dimensional manifolds $(M_a,g_a)$, $a=1,2$. This includes the space-part of the Horowitz-Myers metrics which are conformally compactifiable, have no boundary, and have negative mass, as well as Kottler metrics which are conformally compactifiable and have  totally geodesic boundaries at finite distance.

  Let us denote by $V_a$ the static potential of the asymptotic backgrounds on $M_a$, $a=1,2$.
   For the spherical components of the boundary at infinity, if any, one only needs to consider the static potential  $\sqrt{1+r^2}$,  as follows from \eqref{23VII21.1} below.

   We start with the  boundary-gluing of Section~\ref{s2VII21.1}.
    Letting $x=1/r$ in \eqref{18IX20.4}, and using the fact that the metric is exactly hyperbolic  in $ \mcU_a$ we have
\begin{eqnarray}
 m_a
    &:= &
 m(V_a,\partial M_a)
    =
   - 
   \lim_{x\rightarrow0}\int_{\{x\} \times \partial M_a}   D^j V_a
    ( R{}^i{}_j - \frac {R{}}{n}\delta^i_j)
    \,
    d\sigma_i
    \nonumber
\\
    &= &
   - 
   \lim_{x\rightarrow0}\int_{\{x\} \times
    (\partial M_a\setminus \mcU_a)}   D^j V_a
    ( R{}^i{}_j - \frac {R{}}{n}\delta^i_j)
    \,
    d\sigma_i
     \,.
           \label{22VII21.5}
    \end{eqnarray}
We write $r$ for the radial coordinate as in \eqref{26XI21.1} on $M_1$, so that  the static potential entering into the definition of the mass of $(M_1,g_1)$ reads
\begin{equation}\label{23VII21.1m}
  V_1=  \sqrt{r^2+k}
   \,.
\end{equation}
 Let us write  $(\bar M, \bar g)$ for the boundary-glued manifold,  $\bar V$ for the associated static potential, and $\bar r$ for the radial coordinate  on $\bar M$.
 Viewing $ M_1\setminus \mcU_1$ as a subset of $ \bar M$,  using \eqref{20IX21.11a} with $\psi$ there denoted by $\psi_1$ here,
we can write on $M_1 \setminus \mcU_1$
\begin{equation}\label{23VII21.1}
  \bar V =  {\sqrt{\bar r^2 +\bar k}}
   =   { \bar r } + O(\bar r^{-1})
   =   { \rphi_1  r } + O(  r^{-1})
   =   {V_1}{ \rphi_1    } + O(  r^{-1})
    \,.
\end{equation}
Here $\rphi_1 $ has been extended to the interior of $\bar M_1$ by requiring $\partial_r \rphi_1  =0$.

An identical formula holds on $M_2\setminus \mcU_2$, with a conformal factor which we denote by $\rphi_2$.

Assuming in addition to \eqref{m1}-\eqref{m5} that the $b$-norm of $R{}^i{}_j - \frac {R{}}{n}\delta^i_j$ decays as $o(r^{2-n})$ for large $r$,
  }
  the above leads to
%
\begin{eqnarray}
 \bar m
    &:= &
 m(\bar V ,\partial \bar M ) =
   -
   \lim_{x\rightarrow0}\int_{\{x\} \times \partial \bar M }
      D^j \bar V
    ( R{}^i{}_j - \frac {R{}}{n}\delta^i_j)
    \,
    d\sigma_i
  \nonumber
\\
  &  = &
   - 
   \sum_{a=1}^{2}
   \lim_{x\rightarrow0}\int_{\{x\} \times \partial M_a}
      D^j (\red{\rphi}_a V_a)
    ( R{}^i{}_j - \frac {R{}}{n}\delta^i_j)
    \,
    d\sigma_i
     \,.
           \label{29VII21.1}
    \end{eqnarray}

Comparing with \eqref{22VII21.5}, we see that the mass of the glued manifold is the sum of mass-type integrals, where the original integrands over $\partial M_a$ are adjusted by the conformal factor relating the metric on
the  glued boundary to the original one.

The above applies also to the case where the mass aspect tensor is well-defined, giving a simpler formula
\begin{eqnarray}
 \bar m
  &  = &
   - 
   \sum_{a=1}^{2}
   \int_{ \partial M_a}
      \rphi_a\Big(
  (n-1)\mu_{nn}
+  n \sum_{i=1}^{n-1}\mu_{ii}
   \Big) d^{n-1}\mu_{\hk}
   \,.
           \label{29VII21.1asdf}
    \end{eqnarray}
It should, however, be pointed-out that the gluing-in of an exactly hyperbolic region as in~\cite{ChDelayExotic} is \emph{not known} to lead to a metric with well defined mass-aspect tensor in the region where the deformation of the metric is carried-out, so that in our analysis below we have to use \eqref{29VII21.1}.

 Some comments on the gluing of \cite{ILS} is in order.
It turns out that the difference between the mass formulae for the  boundary-gluing of Section~\ref{s2VII21.1} and that of \cite{ILS} is essentially notational, and concerns the integration domains in \eqref{29VII21.1} and \eqref{29VII21.1asdf}.
  Indeed, after the gluing construction presented in Section 4.2 of \cite{ILS} has been done, the boundary manifold there is obtained by removing  discs $D(\epsilon_a)$ of coordinate radii $\epsilon_a$ from $\partial M_a$, and connecting the boundaries of the discs with a  neck $S^1\times [0,1]$. (The radii $\epsilon_a$ are taken to be one in  \cite[Section~4.2]{ILS}, which is due to a previous rescaling of the coordinates.)
   For the purpose of the formulae below let us  cut the boundary neck $S^1\times [0,1]$ into  two pieces  $N_1:= S^1\times [0,1/2]$ and $N_2:= S^1\times [1/2,1]$, and  let $\Omega_a = (\partial M_a \setminus  D(\epsilon_a) )\cup N_a$.
   In the construction of \cite{ILS} each $\Omega_a$ comes naturally equipped with a constant scalar curvature metric which coincides with the original one on $\partial M_a \setminus  D(\epsilon_a) $.  The relative conformal factors $\psi_a$ are defined on the $\Omega_a$'s with respect to these metrics. We then have
\begin{eqnarray}
 \bar m
    &:= &
 m(\bar V ,\partial \bar M )
  \nonumber
\\
  &  = &
   - 
   \sum_{a=1}^{2}
   \lim_{x\rightarrow0}\int_{\{x\} \times \Omega_a}
      D^j (\red{\rphi}_a V_a)
    ( R{}^i{}_j - \frac {R{}}{n}\delta^i_j)
    \,
    d\sigma_i
     \,.
           \label{29VII21.1ILS}
    \end{eqnarray}
Whenever a  mass aspect tensor is globally defined on the glued manifold, this last formula coincides with
\begin{eqnarray}
 \bar m
  &  = &
   - 
   \sum_{a=1}^{2}
   \int_{ \Omega_a}
      \rphi_a\Big(
  (n-1)\mu_{nn}
+  n \sum_{i=1}^{n-1}\mu_{ii}
   \Big) d^{n-1}\mu_{\hk}
   \,.
           \label{29VII21.1asdfILS}
    \end{eqnarray}

It should be clear how this generalises to the boundary-gluing of several three-dimensional manifolds.

\subsection{Noncompact boundaries}
 \label{ss16IX21}

When defining hyperbolic mass, it is usual to assume, and we did, that the conformal boundary at infinity is compact. This is done for convenience: when the boundary manifold is  compact, to obtain convergence of the integrals defining mass it suffices to impose conditions on the rate of decay of the metric $g$ to the asymptotic background $b$. If the boundary manifold were not compact, one would  need to impose further decay conditions in the non-compact directions at the boundary. However, it is not clear what conditions are relevant for physically interesting fields.

Now, our ``exotic gluings'' construction in~\cite{ChDelayExotic} creates an open neighbourhood, say $\mcO$,  of a subset of the conformal boundary at infinity, such that   the metric is conformal to the hyperbolic one in $\mcO$. One can then create non-trivial non-compact boundary manifolds by removing  closed sets, say $\Sigma$, from the boundary at infinity inside $\partial M\cap \mcO$ and changing the asymptotic background $b$ to a new background $\overline b$ on $\partial M \setminus \Sigma$. A simple
non-trivial example is provided by taking $\partial M$ to be a two-dimensional torus, $\Sigma$ to be a finite collection of points $p_i\in \partial M$ (e.g., one point), and replacing the original background with a flat conformal metric $h_0$ on $\partial M$ by a complete hyperbolic metric with cusps at the (removed) points $p_i$. Then, instead of measuring the mass of $M$ with respect to the original background with a toroidal boundary we can define a mass with respect to the new background $\bar b$, as in \eqref{26XI21.1zd}-\eqref{26XI21.1as}.

Since the removed set $\Sigma$ is contained in the set $\mcO$ where the metric is conformal to the hyperbolic metric, the convergence of the new mass integrals readily follows from the convergence of the original mass: the support of the boundary integrals near infinity is included in a neighborhood of the compact set $\partial M \setminus \overline \mcO$.
The arguments presented in Section~\ref{s5VII21.2} readily lead to the following formula for the new mass:
\begin{eqnarray}
 m(\bar V ,\partial \bar M \setminus \Sigma)  =
   - 
   \lim_{x\rightarrow0}\int_{\{x\} \times (\partial \bar M \setminus \Sigma)}
      D^j (\red{\rphi}  V )
    ( R{}^i{}_j - \frac {R{}}{n}\delta^i_j)
    \,
    d\sigma_i
     \,;
           \label{16XI21.1}
    \end{eqnarray}
recall that
$\psi$ is the conformal factor which relates the original metric $h_k$ to the new one as  $h_{\bar k} = \psi^{-2} h_{k}$,
and $x$ is a coordinate which vanishes at the conformal boundary at infinity.

\section{Higher genus solutions with negative mass}
\label{s29VII21.1}

In this section we  show existence of classes of three-dimensional ALH manifolds $(M,g)$ with constant scalar curvature,  higher-genus conformal boundary at infinity, and negative mass.

\subsection{Genus two}
 \label{ss18VIII21.1}

We start with a construction leading to ALH manifolds with genus at infinity equal to two. Some terminology and preliminary remarks will be needed.

The manifold $(M,g)$ will be obtained by a boundary-gluing of two ALH manifolds, $(M_1,g_1)$ and $(M_2,g_2)$, with  \emph{identical  toroidal conformal geometries at infinity}. We assume existence of a coordinate system near the conformal boundary at infinity in which each of  the original metrics $g_a$ takes the form
\begin{equation}\label{21V22.1}
  g_a = \underbrace{
   \frac{dr^2}{r^2} + r^2 h_0
   }_{\mygb} + e_a
  \,,
  \qquad
   r\ge r_0
   \,,
   \quad
    a=1,2\,,
\end{equation}
for some $r_0>0$, where
%
\begin{equation}\label{21V22.3}
  |e_a|_{\mygb}
  + |D e_a|_{\mygb}
  + |D^2 e_a|_{\mygb}
  \le C r^{-\sigma}
\end{equation}
with constants  $\sigma>5/2$ and $C>0$. This guarantees the possibility of performing an ``exotic gluing'' with well-defined final mass~\cite{ChDelayExotic}, and is satisfied in particular for the Birmingham-Kottler and Horowitz-Myers metrics, in which cases $\sigma=3$.

We make appeal to the construction described in Section~\ref{s2VII21.1}, where  the hyperbolic metric has been glued-in within an $\epsilon$-neighborhood of  boundary points $p_a\in \partial M_a$.
We use the coordinates of \eqref{21V22.1}-\eqref{21V22.3}
 with
polar coordinates for the boundary metric
\begin{equation}\label{29V22/1}
 h_0=d\rho^2 + \rho^2 d\varphi^2
\end{equation}
on $\partial M_a$, with $p_a$ located at the origin of these coordinates.
Such coordinates can always be defined, covering a disc $D(\rho_0)$ for some $\rho_0>0$. Without loss of generality we can rescale the metric $h_0$ on $\partial M_a$ so that the coordinates are defined on the unit disc $D(1)$. (This might rescale the total mass by a positive factor; one then reverts to the original scaling before calculating the mass.)

The metric $g_1$ is thus the original  metric outside the half-ball $\mcU_{1,2\epsilon}$ of coordinate radius $2\epsilon<1$,
and is exactly hyperbolic inside the half-ball $\mcU_{1, \epsilon}$ of coordinate radius $\epsilon$ (see Figure~\ref{F26VIII21.2}); similarly for $g_2$.
 \begin{figure}
	\centering
 \includegraphics{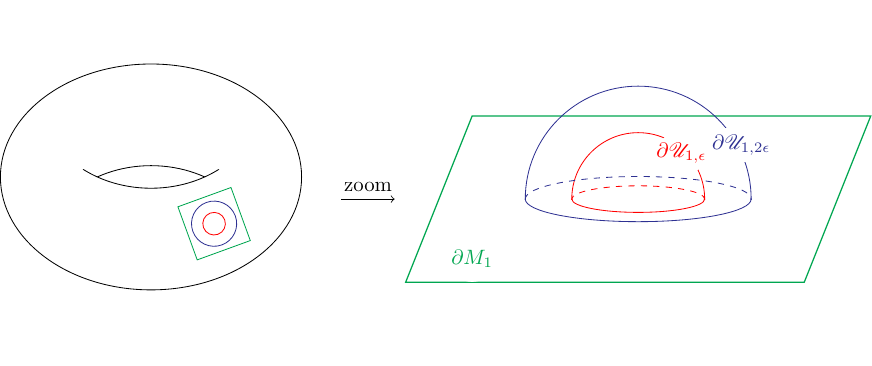}
    \caption{The sets $\mcU_{1, \epsilon}\subset \mcU_{1,2\epsilon}$ and their boundaries.}
	\label{F26VIII21.2}
\end{figure}
For reasons that will become clear in the proof  we need to consider a family of boundary gluings indexed by a parameter $i\to\infty$.
As notation suggests, the parameter $i$ will be taken in $\N$; we will indicate how to transition to a continuous parameter $i$ in Remark~\ref{R25V22.1} below.
For $i\ge \red{8}/\epsilon$  we choose both hyperbolic hyperplanes $\blue{\horo_{1,i}}$ and $\horo_{2,i}$ of Section~\ref{s2VII21.1}
to be  half-spheres of radius $1/i$ centered at the origin of the coordinates \eqref{29V22/1}.
 We choose any pair $(\Lambda_1:=\Lambda_{\horo_{1,i},+},\Lambda_2:=\Lambda_{\horo_{2,i},-})$ as in Section~\ref{s2VII21.1} to obtain the boundary-glued manifold $M:=M_{\Lambda_1,\Lambda_2}$.

Recall that, given a flat torus $(\T^2,h_0)$ and a point $p \in \T^2$, there exists on $\T^2\setminus \{p\}$ a smooth function $\ominf $ such that the metric $e^{\ominf }h_0$ is complete, has constant Gauss curvature equal to $-1$, with $(\T^2,e^{\ominf }h_0)$  having finite total area; see Figure~\ref{F21XI21.1} and Appendix~\ref{sApp20XI21.1}, compare~\cite[Proposition~2.3]{GuillarmouSurfaces}.
This metric will be referred to as the \emph{hyperbolic metric with a cusp at $p$.}
 \begin{figure}
	\centering
 \includegraphics{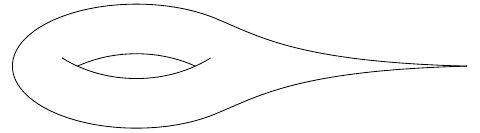}
    \caption{A punctured torus with a hyperbolic metric; the figure is misleading in that the cusp region is actually infinitely long.}
	\label{F21XI21.1}
\end{figure}

We claim:

\begin{Theorem}
\label{T29VII21.1}
Consider two three-dimensional ALH manifolds $(M_a,g_a)$, $a=1,2$, with constant scalar curvature  and with a metric of the form \eqref{21V22.1}-\eqref{21V22.3} with the same flat metric $h_0$.
Let $p\in \partial M_1$ and let $e^{\ominf } h_0$ be the unique constant negative curvature metric with a cusp at $p$. The mass of the Maskit-doubled metric as described above converges, as $\epsilon$ tends to zero and $i$ tends to infinity, to the finite limit
\begin{equation}\label{23V22.1}
 - 
   \sum_{a=1}^{2}
   \lim_{r\rightarrow \infty}\int_{\{r\} \times \T^2}
      D^j (e^{-\ominf/2}r )
    ( R{}^i{}_j - \frac {R{}}{n}\delta^i_j)
    \,
    d\sigma_i
   \,.
\end{equation}
In this equation $R_{ij}$ denotes the Ricci tensor of $g_1$ when $a=1$ and of $g_2$ when $a=2$.
\end{Theorem}

This leads us to existence of negative-mass solutions without boundary:

\begin{Corollary}
\label{C21V22.1}
There exist three-dimensional conformally compactifiable ALH manifolds without boundary at finite distance, with constant scalar curvature, genus two infinity, and negative mass.
\end{Corollary}

\noindent{\sc Proof of Corollary~\ref{C21V22.1}:}
We apply Theorem~\ref{T29VIII21.1} to two identical copies of the space-part of a Horowitz-Myers metric. Keeping in mind that $\mc>0$,  it follows from \eqref{23V22.1} and \eqref{29VII21.1asdf}  that $m$ approaches
\begin{equation}\label{23V22.111}
 -\infty< - 2\mc \int_{\T^2} e^{-\ominf /2} d\mu_{h_0}
   < 0
\end{equation}
when $\epsilon$ tends to zero and $i$ tends to infinity, and hence is negative for sufficiently small $\epsilon$ and sufficiently large $i$.
\qed

\begin{Remark}
  \label{R23V22.1}
 {\rm
 In the case of the square torus with area equal one, discussed in Appendix~\ref{sApp20XI21.1}, the function $\ominf$ in \eqref{23V22.111} satisfies (cf.~\eqref{20XI21.9})
 \begin{equation}\label{20XI21.8}
   e^{-\ominf /2} \le  \frac{\Gamma(1/4)^2}{4 \pi ^{3/2}} \approx 0.59
   \,,
 \end{equation}
 which leads to a rough lower estimate of the mass $m$ of the manifolds obtained by gluing two identical copies of the space-part of the Horowitz-Myers  metric, for small $\epsilon$ and large $i$,
 \begin{equation}\label{20XI21.7}
 -1.18\, \mc<  m  <0
   \,.
 \end{equation}
 %
 An upper bound  can be obtained by numerically integrating the inverse of the right-hand side  of \eqref{20XI21.11}, leading to
 \begin{equation}\label{20XI21.71}
  -1.18 \,\mc
  <  m < - 0.68 \, \mc
   \,.
 \end{equation}
 (The   integral of the inverse of the left-hand side  of \eqref{20XI21.11} gives a worse lower bound $\approx -4.02 \,\mc$.)
 }
 \end{Remark}

We also obtain negative-mass solutions with a toroidal black-hole boundary:

\begin{Corollary}
\label{C21V22.1a}
There exist three-dimensional conformally compactifiable ALH manifolds with a connected toroidal minimal boundary, with constant scalar curvature, genus two infinity, and negative mass.
\end{Corollary}

\noindent{\sc Proof of Corollary~\ref{C21V22.1a}:}
Apply Theorem~\ref{T29VIII21.1} to a Maskit gluing of a Horowitz-Myers metric with mass parameter $\mc{}_{,1}>0$ and with a Birmingham-Kottler metric with a toroidal black hole with mass parameter $0< \mc{}_{,2}<\mc{}_{,1}/2$,    so that the limiting total mass (see \eqref{29VII21.1asdf}, cf.~\eqref{3VIII21.2}  and \eqref{25V22.98})
$$
  (2 \mc{}_{,2}-\mc{}_{,1}) \int_{\T^2} e^{-\ominf /2} d\mu_{h_0}
$$
 is negative.
\qedskip

\noindent{\sc Proof of Theorem~\ref{T29VII21.1}:
}
By construction the boundary $\partial M$  of the new manifold is a doubling of
$$
  \zMone := \partial M_1\setminus D(\red{1/i})
   \,.
$$
We will often view both
 $\zMone $ and its double as   subsets of $\partial M$.

Consider the conformal class of metrics on $\partial M$ induced by $g$.
This conformal class depends upon $i$ but is independent of $\epsilon$, keeping in mind that $\epsilon \ge \red{8}/i$.
In this class there exists  a unique metric $h_{-1}$ (depending upon $i$) with constant scalar curvature equal to minus two which can be constructed as follows:
Let $\hath_i   $
be any representative of the conformal class of metrics on $\partial M_1$ induced by $g$
  which
  is

 \begin{enumerate}
   \item
   invariant under reflection across $\partial\zMone $, so that  the coordinate circle
$$
 S_{\red{1/i}}=\partial\zMone
$$
of coordinate radius $\red{1/i}$  is a totally geodesic boundary for $(\zMone,\red{\hath_i})$,  and
  which
  is
   \item  invariant under rotations near $S_{\red{1/i}}$.
 \end{enumerate}

 While any metric $\red{\hath_i}$ as in 1.\ and 2.\ is adequate, a useful
 example is provided by the metric which on $\partial M_1$ equals $\chi^2 h_0$, where $\chi\ge 1/2$ is any smooth function on $\partial M_1\setminus \{p\}$ which equals to $1/\rho$ on $D(1)\setminus \{0\}$,
where $\rho$ is the coordinate radius in the local coordinates   on $\partial M_1$. This can be accompanied by the introduction of a new coordinate
\begin{equation}\label{25XI21.1}
  \hrho = \frac{ \red{\log} ( \rho)}{ \red{\log} (i )}+1
   \,,
\end{equation}
so that the flat metric
$$
 h_0 = d\rho^2 + \rho^2 d\varphi^2
 \,,
 \qquad
  \rho\in [1/i,1]
$$
becomes
\begin{equation}\label{25XI21.2}
 h_0 =\rho^2
  \big( \frac{d\rho^2}{\rho^2}+ d\varphi^2\big)
  = \rho^2
  \big ( {\log}^2(i) d\hrho^2 +  d\varphi^2 \big)
 \,,
 \qquad
  \hrho\in [0,1]
  \,,
\end{equation}
thus $\red{\hath_i} = \red{\log}^2(i) d\hrho^2 + d\varphi^2$ on $[1,i]\times S^1$.

In this coordinate system the \emph{reflection across $\partial\zMone $} is the map $(\hrho,\varphi)\mapsto (-\hrho,\varphi)$, after extending the range of $\hrho$ from $[0,1]$ to $[-1,1]$.

Given a function $f$  defined on $\overline{D(1)\setminus D(1/i)}$ in the original coordinates $(\rho,\varphi)$, the \emph{mirror-symmetric doubling} of $f$ is defined as the function
\begin{equation}\label{2XI21.3}
  [-1,1]\times S^1 \ni (\hrho,\varphi)\mapsto
\hat f (\hrho,\varphi)=
  \left\{
    \begin{array}{ll}
      f( e^{ \red{\log}(i)(\hrho -1) },\varphi) & \hbox{$\hrho \ge0$;} \\
      f(e^{\red{\log}(i)(-\hrho  -1) },\varphi) & \hbox{$\hrho <0$.}
    \end{array}
  \right.
\end{equation}
The mirror-symmetric doubling of a tensor field is obtained by transforming the tensor field to the coordinates $(\hrho,\varphi)$  and mirror-doubling all its coordinate-components. For example, after doubling of $\zMone$, the tensor field  $\red{\hath_i} = \red{\log}^2(i) d\hrho^2 + d\varphi^2$ just defined above on $[0,1]\times S^1$ maintains the same form $ \red{\log}^2(i) d\hrho^2 + d\varphi^2$ on $[-1,1]\times S^1$ after doubling.

Another example, which plays an important role in our proof below, is provided by the daunting-looking metric (cf., e.g., \cite{MazzeoTaylor})
\begin{equation}
 \label{9XI21.3a}
  e^{\mv_{*,i}}h_0  :=
\Big(
 \frac{\pi }{
 \log (i^2) \sin( \pi \frac{\log \rho }{\log  (i^2)} )  \rho
}
 \Big)^2 (d\rho^2+\rho^2 d\varphi^2)
  \,, \qquad
  \rho\in (1/i^2,1 )
  \,,
\end{equation}
which has constant negative  scalar curvature equal to $-2$.
The conformal factor $e^{\mv_{*,i}}$ tends to infinity  at $\rho=\rhozero /i^2$ and at $\rho=\rhozero $. In the coordinates $(\hat \rho, \varphi)$ the metric \eqref{9XI21.3a} reads
\begin{equation}
 \label{25XI21.6}
  e^{\mv_{*,i}}h_0  =
\Big(
 \frac{\pi }{
 2  \cos( \pi \hrho/2  )
}
 \Big)^2 (d\hrho^2+ \frac{1}{\red{\log^2(i)}} d\varphi^2)
  \,, \qquad
  \hrho\in (-1,1)
  \,,
\end{equation}
which is manifestly mirror-invariant. Note that the circle $\hrho =0$ is a closed geodesic minimising length, of length $\pi^2/\red{\log (i)}$.

Let $u_i:\partial M\to \R $ be the unique solution of the two-dimensional Yamabe equation
%
\begin{equation}\label{31VII21.31}
   \Delta_{\red{\hath_i}  }u_i  =   - R e^{ u_i} +   \hatR_i
  \,,
\end{equation}
with $R=-2$,  and where $\hatR_i  $ is the scalar curvature of the metric $\red{\hath_i}  $, so that the metric $e^{ u_i}\red{\hath_i}  $ has scalar curvature $R$.
It is important in what follows that the function  $u_i$ is independent of $\epsilon$.

(As uniqueness is likewise important in our analysis, let us provide an argument: let $\tilde u_i$ and $\hat u_i$ be two solutions, rename the metric $ e^{\hat u_i}\red{\hath_i}$ to $\red{\hath_i}$, then from $\tilde u_i$ one obtains a solution of \eqref{31VII21.31} with
 $R=\hatR=-2$. Multiplying the resulting equation by $u_i$ and integrating over $\partial M$ one obtains
\begin{equation}\label{31VII21.35}
   \int |du_i|^2 d\mu_{\red{\hath_i}}  =   2 \int (1- e^{u_i})u_id\mu_{\red{\hath_i}}
  \,.
\end{equation}
Since the integrand of the right-hand side is negative except at $u_i=0$, we find $u_i\equiv 0$, hence $\tilde u_i\equiv \hat u_i$.)

Uniqueness  implies that $u_i$ is invariant under reflection across $\zMone $.
Hence $u_i$ has  vanishing normal derivative on $S_{\red{1/i}}$.

As $\red{\hath_i}  $ is conformal to $h_0$ on $\zMone$, there exists a
 function $\red{\hat u_i} $
so that we have
$$
 \red{\hath_i}   =e^{   \red{\hat u_i} } h_0
$$
on $\zMone$.
Then the metric
$$
 e^{ u_i+\red{\hat u_i} } h_0
$$
defined on $\zMone$, has scalar curvature equal to minus two. The mirror-symmetric doubling of $\red{\hath_i}$ across $S_{\red{1/i}}$  coincides with the metric $h_{-1}$ on $\partial M_1$.

Now, by construction, $\red{\hat u_i} $ is rotation-invariant near $S_{\red{1/i}}$. The fact that $S_{\red{1/i}}$ is totally geodesic  is equivalent to the  vanishing of the radial derivative of
$\red{e^{ \red{\hat u_i} /2}\rho}$  on $S_{\red{1/i}}$:
\begin{equation}\label{1VIII31.1}
 0 =  \partial_\rho(\red{e^{ \red{\hat u_i} /2}\rho})|_{S_{\red{1/i}}}
 \qquad
  \Longrightarrow
 \qquad
 \partial_\rho \red{\hat u_i}  |_{S_{\red{1/i}}} = \red{-2i}
  \,.
\end{equation}
From what has been said we  have
\begin{equation}\label{1VIII31.2}
 \partial_\rho   u_i |_{S_{\red{1/i}}} = 0
  \,.
\end{equation}
Since $h_0$ is flat, we conclude that the function
$$
 \mv_i:=  u_i+\red{\hat u_i}
$$
satisfies on $\zMone$ the equation
\begin{equation}\label{31VII21.32}
   \Delta_{h_0  }\mv_i  =   2 e^{\mv_i}
  \,,
\end{equation}
with Neumann boundary data at the coordinate circle $S_{\red{1/i}}$ of radius $\red{1/i}$ centered at the origin:
\begin{equation}\label{1VIII31.3}
 \partial_\rho   \mv_i |_{S_{\red{1/i}}} = -\red{2i}
  \,.
\end{equation}

Now,
while the function $\mv_i$ depends  only upon $i$, the metric on $M$ depends upon both   $\epsilon$ and $i$. This has the unfortunate consequence  that the sign of the mass of $M$ is not clear. In what follows we will determine this sign for $i$ large and $\epsilon$ small. For this we need to understand what happens with the mass integrand \eqref{29VII21.1} when $i$ tends to infinity and $\epsilon$ tends to zero.
 The idea of the argument is to show that the sequence of functions $\mv_i$, or at least a subsequence thereof, converges to a limit,  with sufficient control of the limit to guarantee control of the mass. The needed convergence result is perhaps contained in \cite{Wolpert}, but it is not completely apparent to us that this is the case, so we provide a direct  argument, different from the one in~\cite{Wolpert}.

The maximum principle shows that $\mv_i$ has no interior \mage{maximum} on  the compact manifold with boundary $ \partial M_1 \setminus D(a)$ for $a\in[1/i,1]$, where $D(a)$ denotes an open coordinate disc of radius $a$  lying on the boundary $\partial M_1$ and centered at $p$.

Integrating \eqref{31VII21.32} over $\zMone$ and using \eqref{1VIII31.3} one finds
%
\begin{equation}\label{1VIII31.4}
  \int_{\zMone} e^{\mv_i} d\mu_{h_0}= \red{2\pi}\,.
\end{equation}
%

We continue the proof of Theorem~\ref{T29VII21.1} with the (well known, cf.\ e.g.~\cite{MazzeoTaylor}) observation, that solutions of the equation
\begin{equation}\label{23X21.41}
 \Delta \mv= 2 e^{\mv}
\end{equation}
satisfy a comparison principle, which we formulate in the simplest form sufficient for our purposes: If both $\mv$ and $\hat \mv$ satisfy \eqref{23X21.41} in a bounded domain $\Omega$ with  boundary, and if $\hat \mv >\mv$ near the boundary, then  $\hat \mv > \mv$ on $\Omega$. (This comparison is particularly useful with functions $\hat \mv$ which tend to infinity when $\partial\Omega$  is approached.)
Indeed, if $\omega$ is continuous on the closure $\bar \Omega$ of $\Omega$, the function $\mv-\hat \mv$ is \mage{then negative} near $\partial\Omega$, and satisfies the equation
\begin{equation}\label{14IX21.1a}
   \Delta_{h_0  } (\mv -\hat \mv )  =    2 e^{ \mv  }-2 e^{\hat \mv }
  \,.
\end{equation}
If  $ \mv > \hat  \mv $ somewhere, then $\mv-\hat \mv$  would have a positive maximum at some point $q\in \Omega $  away from the boundary. Since the function \mage{$x\mapsto 2 e^{x}$ }is increasing, the  right-hand side of \eqref{14IX21.1a} would be positive at $q$. But the left-hand side is nonpositive at a maximum, a contradiction.

We rephrase this loosely as
%
\begin{equation}\label{14IX21.2as}
  \hat \mv    > \mv \ \mbox{on $\partial\Omega$}
   \quad
   \Longrightarrow
   \quad
     \hat \mv   > \mv\ \mbox{on $\Omega$.}
\end{equation}
%



For $0<a<b$ let
$$
\Gamma(a,b) := {D(b) \setminus D(a)}
\,.
$$
The comparison principle allows us to prove:

\begin{Lemma}
\label{L23X21.1m}
On $\Gamma(1/i,\rhozero )$ it holds that
\begin{equation}
 \label{9XI21.2}
e^{w_i}\leq  e^{\mv_{*,i}}:=
\Big(
 \frac{\pi }{
 \log (i^2) \sin( \pi \frac{\log \rho }{\log  (i^2)} )  \rho
}
 \Big)^2
\,.
\end{equation}
%
%
\end{Lemma}

\begin{remark}
 \label{R22V22.1}
{\rm
On the circle $\rho=1/i$ we have
$$
e^{\mv_{*,i}}=\frac{\pi^2 i^2}{4\log^2(i)},
$$
which is unbounded in $i$, but the metric length $\ell_i$ of $S_{1/i}$ equals
$$
\ell_i=\frac1i\int_{\varphi\in[0,2\pi]}(e^{\mv_{i}/2})_{|\rho=1/i}d\varphi\leq\frac1i\int_{\varphi\in[0,2\pi]}(e^{\mv_{*,i}/2})_{|\rho=1/i}d\varphi=\frac{\pi^2}{\log{i}},
$$
so that $\ell_i$ approaches zero as $i$ tends to infinity.
}
\end{remark}

{\noindent \sc Proof of Lemma~\ref{L23X21.1m}:}
The metric  (compare \eqref{9XI21.3a})
\begin{equation}
 \label{9XI21.3}
 e^{\mv_{*,\sqrt t}} h_0  = \big(
  \frac{\pi}{|\log t|} \csc( \pi \frac{\log |z|)}{\log |t|} \big)^2 \frac{|dz|^2}{|z|^2}
\end{equation}
%
 has constant negative  scalar curvature equal to $-2$. Furthermore, the circle $S_{1/i}$ is minimal when $t=i{^2}$:
\begin{equation}\label{18XI21.1AA}
  \partial_\rho (\rho e^{\frac{\mv_{*,i}}2})|_{\rho = 1/i} = 0
  \,.
\end{equation}
The conformal factor $e^{\mv_{*,i}}$ tends to infinity  at $\rho=\rhozero /i^2$ and at $\rho=\rhozero $.

The change of the complex variable $z\mapsto   w=  \rhozerosquared /(i^2  z)$
reproduces the metric on $w\in \Gamma(\rhozero /i^2,\rhozero )$ and exchanges the regions on both sides of the geodesic $|z|=|w|=\rhozero /i$.
%
Our  doubled  metric on $\partial M$  is likewise symmetric  relatively to this geodesic. The  conformal factor $e^{\mv_i}$, extended to $\rho\in (1/i^2,1/i)$ using the map just described,  provides a smooth solution of the Yamabe equation \eqref{23X21.41} on $\Gamma(1/i^2,1)$ (compare the discussion around \eqref{2XI21.3}). It takes finite values both at $\rho=\rhozero /i^2$ and at $\rho=\rhozero $. The result follows from the comparison principle.
\qed

\begin{Corollary}
 \label{C21XI21} For any $\rho_1\in(0,1)$,
 there exists a constant $\hat c =\hat c(\rho_1)$  such that
 $$\mv_i \le \hat  c
 $$
 on $\zMone \setminus D(\rho_1)$, independently of $i$.
\end{Corollary}

\proof
The inequality
\begin{equation}
 \label{9XI21.2b}
e^{w_i}\leq  e^{\mv_{*,i}}=
\Big(
 \frac{\pi }{
 \log (i^2) \sin( \pi \frac{\log \rho_1 }{\log  (i^2)} )  \rho_1
}
 \Big)^2\longrightarrow_{i\rightarrow +\infty}\frac{1}{\rho_1^2\log^2(\rho_1)}
\,,
\end{equation}
shows that the $\mv_i$'s are bounded by a constant $\hat c(\rho_1 )$ independently of $i$ on $S(\rho_1)$ for $\rho_1\in (0,1)$.
The result follows from the maximum principle.
\qed

The corollary gives an estimation of the conformal factors from above. In order to prove convergence of the sequence $\mv_i$, away from the puncture, we also need to bound the sequence of conformal factors away from zero. The next lemma provides a first step towards this. At its heart lies the
inequality \eqref{6X21.2} which, together with \eqref{18XI21.1},
shows that the area does not concentrate near the minimal geodesic $S_{1/i}$.
The parameter $\repsilon>0$ in Lemma~\ref{L23X21.1} below should not be confused with the parameter $\epsilon$ introduced by the exotic gluing of $M_1$ with $M_2$:

\begin{Lemma}
\label{L23X21.1}
For all $\repsilon>0$ sufficiently small and all $i$ sufficiently large we have the bounds
\begin{equation}\label{15VIII21.21b}
  \inf_{\red{\partial M_1  \setminus  D(\repsilon/2)}} (\mv_i - \red{\log} 2)\le
   \red{\log}  \frac{\red{\pi}}{|\partial M_1 \setminus \red{D(}\repsilon/2)|_{h_0}}
   \le \sup_{\red{\partial M_1  \setminus  D(\repsilon/2)}} \mv_i
   \,,
\end{equation}
where $|U|_{h_0}$ denotes the $h_0$ area of a set $U$,
and note that the middle term is independent of $i$.
\end{Lemma}
\proof
We have
\begin{equation}\label{18XI21.1A}
  \int_{\Gamma(a,b)} e^{\mv_{*,i} }  d\mu_{h_0}
  = -\frac{2 \pi ^2 \cot \left(\frac{\pi
    \log (\rho )}{\log
   \left(i^2\right)}\right)}{\log
   \left(i^2\right)}\bigg|^b_a
   \,,
\end{equation}
with
\begin{equation}\label{18XI21.1}
  \int_{\Gamma(1/i,\repsilon/2)} e^{\mv_{*,i} }  d\mu_{h_0}
  = \frac{2 \pi ^2 \cot \left(\frac{\pi
   \log (2/\repsilon)}{\log
   \left(i^2\right)}\right)}{\log
   \left(i^2\right)}
   \to_{i\to\infty}
   \frac{2\pi }{\log (2/\repsilon)}
   \to_{\repsilon\to 0} 0
   \,.
\end{equation}
Equation \eqref{1VIII31.4} gives
\begin{eqnarray}
 2\pi
  &= &
  \int_{\zMone} e^{\mv_i} d\mu_{h_0}=    \int_{\red{\partial M_1  \setminus  D(\repsilon/2)}}  e^{\mv_i} d\mu_{h_0}
   +
   \int_{\Gamma(1/i,\repsilon/2) } e^{\mv_i} d\mu_{h_0}
  \label{6X21.2ax}
   \,.
  \label{6X21.21}
\end{eqnarray}
The estimate \eqref{9XI21.2} shows that
%
\begin{eqnarray}
\int_{\Gamma(1/i,\repsilon/2) } e^{\mv_i} d\mu_{h_0}
  &\le&
     \int_{\Gamma(1/i,\repsilon/2) } e^{\mv_{*,i} } d\mu_{h_0}
   \,.
  \label{6X21.2}
\end{eqnarray}
It follows from \eqref{18XI21.1} that there exists $\repsilon_0$ such that for all   $2/i <\repsilon\le \repsilon_0$ the last term in \eqref{6X21.2ax} is in $(0,\pi)$, 
which implies
\begin{eqnarray}
     \pi <  \int_{\red{\partial M_1  \setminus  D(\repsilon/2)}}  e^{\mv_i} d\mu_{h_0}
 < 2 \pi
   \,.
  \label{6X21.3}
\end{eqnarray}
The conclusion readily follows using
\begin{eqnarray}
      \int_{\red{\partial M_1  \setminus  D(\repsilon/2)}}  e^{\inf \mv_i } d\mu_{h_0}
 \le  \int_{\red{\partial M_1  \setminus  D(\repsilon/2)}}  e^{\mv_i} d\mu_{h_0}
 \le  \int_{\red{\partial M_1  \setminus  D(\repsilon/2)}}  e^{\sup \mv_i}  d\mu_{h_0}
   \,.
  \label{6X21.367}
\end{eqnarray}
\qedskip

We continue with:

\begin{Lemma} \label{L23X21.2}
There exists a smooth function
$$
 \mv_\infty:\partial M_1\setminus \{p\} \to \R
$$
such that a subsequence of  $\{\mv_{i_j}\}_{j\in\N}$ converges uniformly, together with any number of derivatives, to $\mv_\infty$ on every compact subset of $\partial M_1\setminus \{p\}$.
\end{Lemma}

\proof
Let $K$ be any compact subset of $\partial M_1\setminus \{p\}$, there exists $\rho_K>0$ such that $K\subset \partial M_1\setminus D(\rho_K)$. It thus suffices to prove the result with $ K=\partial M_1\setminus D(\rho_K)$, which will be assumed from now on.

Let $K_1 =  \partial M_1\setminus D(\rho_K/2)$.
Choosing
$\repsilon<\rho_K /4$
in Lemma~\ref{L23X21.1} ensures that  \eqref{15VIII21.21b}  holds  on $K_1$ for all   $i\ge i_1 $ for some $i_1<\infty$.

Let $i\ge i_1$, by  Corollary~\ref{C21XI21}
there exists a constant $c_1$, independent of $i$, such that
\begin{equation}\label{24X21.1}
 \red{v_i}:=-\mv_i \ge c_1
\end{equation}
on $K_1$. Define
$$
 \red{\hat v_i} : = \red{v_i} -c_1 +1
 \,.
$$
It holds that $\red{\hat v_i} \ge 1$ on $K_1$. Moreover, $\red{\hat v_i}$                                                                     satisfies the equation%
\footnote{We are grateful to   Yanyan Li for providing the  argument leading to \eqref{6VIII21.11a}.}
\begin{equation}\label{VIVIII21.1}
  \Delta_{h_0} \red{\hat v_i} = \red{\psi_i}\red{\hat v_i}
  \,,
\end{equation}
where
\begin{equation}\label{VIVIII21.2}
0 \ge \red{\psi_i}  : =-  2 \frac{e^{- \red{v_i}}}{  \red{\hat v_i}} \ge -2  e^{-c_1}
  \,.
\end{equation}
By Harnack's inequality there exists a constant $C_1= C_1(K,K_1) >0$ such that on $K$ we have
\begin{equation}\label{6VIII21.11m}
  \sup_K \red{\hat v_i} \le C_1 \inf_{K_1} \red{\hat v_i}
  \,.
\end{equation}
This, together with the definition of $\red{\hat v_i}$, shows that
\begin{equation}\label{6VIII21.11a}
  \sup_K \red{v_i} \le C_1 \inf_{K_1} \red{v_i} +d_1
  =
  C_1 (-\sup_{K_1} {\mv_i}) +d_1
  \,,
\end{equation}
for some constant $d_1$.
Equation~\eqref{15VIII21.21b} gives
\begin{equation}\label{15VIII21.21bxdv}
-  \sup_{K_1} \mv_i
    \le - \ln  \frac{\red{\pi}}{|\partial M_1 \setminus \red{D(}\repsilon/2)|_{h_0}}
    =:c_2
   \,,
\end{equation}
From \eqref{6VIII21.11a} we obtain
\begin{equation}\label{23X21.51}
   \sup_K \red{v_{i }} \le   C_1  c_2 +d_1
   \quad
   \Longrightarrow
   \quad
     \inf_K \red{\mv_{i }} =  - \sup_K \red{v_{i }} \ge  -(C_1  c_2 +d_1)
   \,.
\end{equation}
This, together with \eqref{24X21.1}, shows that
that for every compact subset $K$ of $\partial M_1\setminus \{p\}$
there exists a constant $\hat C_K$ such that
\begin{equation}
\label{9VIII211}
 - \hat C_K \le  \red{\mv_{i }} \le \hat C_K
 \,.
\end{equation}
Elliptic estimates,  together with a standard diagonalisation argument, show   that there exists a subsequence $\mv_{i_{j }}$  which
converges uniformly on every compact subset of $\partial M_1 \setminus \{p\}$ to a solution $\mv_\infty $
of  \eqref{31VII21.32} on $\partial M_1 \setminus \{p\}$. Convergence of derivatives follows again from elliptic estimates.
\qedskip

We return to  the proof of  Theorem~\ref{T29VII21.1}.
Equation~\eqref{1VIII31.4} shows that we can use  \cite{ChouWan2} to conclude that there exists
 $\alpha >-2 $ such that for small $\rho$
we have either
\begin{equation}\label{2VIII31.1}
  \mv_\infty  = \alpha  \red{\log}  \rho + O(1)
  \qquad
   \Longrightarrow
  \qquad
  e^{\mv_\infty } \sim \rho^{\alpha}
  \,,
\end{equation}
or
\begin{equation}\label{2VIII31.3}
  \mv_\infty  = -2  \red{\log}  (- \rho \red{\log}  \rho ) + O(1)
  \qquad
   \Longrightarrow
  \qquad
  e^{\mv_\infty } \sim \rho^{-2} \red{\log} ^{-2} \rho
  \,.
\end{equation}
%
Now, we claim that the case $\alpha\ge 0$ cannot occur in our context. Indeed, if \eqref{2VIII31.1} holds, integrating \eqref{23X21.41} over
$\T^2\setminus D(p,\delta)$ one finds, for $i$ large enough,
\begin{equation}\label{23X21.41a4}
  0<  2 \int_{\T^2\setminus D(p,\delta)} e^{\mv_i}
   \,
     d\mu_{h_0}
 = \int_{\T^2\setminus D(p,\delta)} \Delta \mv_i
   \,
   d\mu_{h_0}
 = - \oint_{S_{\delta}}  \partial_\rho \mv_i \, \rho
   \,
   d\varphi\to_{i\to\infty} -2\pi\alpha + O(\delta)
 \,.
\end{equation}
Choosing $\delta $ sufficiently small, we conclude that
\begin{equation}\label{22V221}
  \alpha \in (-2,0)
   \,.
\end{equation}

We are ready now to analyse the mass  of $M$. By construction, $M$ is the union of  $M_1\setminus \mcU_{1,\red{1/i}}$ and  $M_2\setminus \mcU_{1,\red{1/i}}$, identified along the boundary $\blue{\horo_{1,i}}$.
 Since the metric is exactly hyperbolic in $\mcU_{1,\epsilon}$, the mass integrands vanish  in $\mcU_{1,\epsilon}\setminus \mcU_{1,1/i}$.
  Let $m_{i,\epsilon}$ denote the mass of either summand in \eqref{29VII21.1}.
In space-dimension $n=3$ we have
\begin{eqnarray}
 \red{m_{i,\epsilon}}
  &  = &
   -
   \lim_{x\rightarrow0}\int_{\{x\} \times \partial M_1}
   D^j ( \rphi_i V_1 )
    ( R{}^k{}_j - \frac {R{}}{3}\delta^k_j)
    \,
    d\sigma_k
     \nonumber
\\
  &  = &
   -
   \lim_{x\rightarrow0}\int_{\{x\} \times
   (\red{\partial M_1\setminus D(\epsilon)} ) }
   D^j (\rphi_i V_1)
    ( R{}^k{}_j - \frac {R{}}{3}\delta^k_j)
    \,
    d\sigma_k
     \,,
           \label{2VIII21.2mk}
    \end{eqnarray}
where
\begin{equation}\label{19IX21.1ac}
  \rphi_i = e^{-\mv_i/2}
  \,.
\end{equation}

We will need  to estimate  the derivatives of the functions   $\psi_i$. As a step towards this,
for $y\in  \overline{D(4)\setminus D(1/2)}$ we set
\begin{equation}\label{15VIII21.12}
 \neww (y) =
 \left\{
   \begin{array}{ll}
 \red{\mv_i} (\epsilon  y) - \alpha \red{\log}  (\epsilon  |y|), & \hbox{under \eqref{2VIII31.1};} \\
    \red{\mv_i} (\epsilon  y) +2  \red{\log}  (- \epsilon  |y|\red{\log}  (\epsilon  |y|)), & \hbox{under \eqref{2VIII31.3}.}
   \end{array}
 \right.
\end{equation}
Each of the functions $\neww $  satisfies on $D(4)\setminus \overline{D(1/2)}$ the equation
\begin{equation}\label{15VIII21.14}
  \Delta_{h_0  }\neww =
 \left\{
   \begin{array}{ll}
   2 \epsilon^{2+\alpha } |y|^{\alpha}    e^{\neww}, & \hbox{under \eqref{2VIII31.1};}
   \\
 \displaystyle
     \frac{2  e^{\neww} }{ |y|^2(\red{\log}  (\epsilon  |y|))^2}-\frac{2}{|y|^2(\red{\log}  (\epsilon  |y|))^2}
 , & \hbox{under \eqref{2VIII31.3}.}
   \end{array}
 \right.
\end{equation}
For every $\epsilon$ there exists $i_0(\epsilon)$ such that  the functions $\neww$ are bounded uniformly in $\epsilon$ and $i$ for $i \ge i_0$.
This leads to an estimate on the derivatives of $\neww$ on $\overline{D(2)\setminus D(1)}$ (see~\cite[Section~3.4]{GT}):
\begin{equation}\label{22VII21.11}
  |\partial\neww |+ |\partial \partial\neww |
   \le C
   \,,
\end{equation}
for some constant $C$ independent of $\epsilon$ and of $i$ provided that $i\ge i_0(\epsilon)$.

In view of \eqref{2VIII21.2mk}-\eqref{19IX21.1ac}, we will need estimates for
$e^{-\mv_i/2}$.
From  \eqref{22VII21.11} one obtains
on $\overline{D(2\epsilon)\setminus D(\epsilon)}$, for  $i\ge i_0(\epsilon)$,
\begin{equation}\label{15VIII31.2}
  | \partial _A
  e^{ \red{-\mv_i/2}} | \le C \epsilon^{-1} e^{ \red{-\mv_i/2}}
  \,,
  \qquad
    | \partial _A \partial_B
  e^{ \red{-\mv_i/2}} | \le  C\epsilon^{-2} e^{ \red{-\mv_i/2}}
  \,,
\end{equation}
for some possibly different constant $C$ which is independent of $\epsilon$ and $i$ when $i\ge i(\epsilon)$.

Directly from~\cite{ChouWan2},
or replacing $v_i$ by $v_\infty$  in the argument just given,  one finds
\begin{equation}\label{25VIII21.1}
  | \partial _A
  e^{ \red{-\mv_\infty/2}} | \le C \rho^{-1} e^{ \red{-\mv_\infty/2}}
  \,,
  \qquad
    | \partial _A \partial_B
  e^{ \red{-\mv_\infty/2}} | \le C \rho^{-2} e^{ \red{-\mv_\infty/2}}
  \,,
\end{equation}
for possibly yet another constant $C$.

Equation~\eqref{2VIII21.2mk} can be rewritten as
\begin{eqnarray}
 \red{m_{i,\epsilon}}
  &  = &
   \underbrace{
    - 2  \int_{ \red{\partial M_1\setminus D(2\epsilon)}}
    D^j ( \rphi_i V_1)
    ( R{}^k{}_j - \frac {R{}}{n}\delta^k_j)
    \,
    d\sigma_k
    }_{=: \hat m_{i,\epsilon}}
    \nonumber
\\
 &&
      - 2
   \underbrace{
    \lim_{x\rightarrow0}\int_{\{x\} \times
    (\blue{\overline{D(2\epsilon)\setminus D(\epsilon/2)}})}
    D^j ( \rphi_i V_1)
    ( R{}^k{}_j - \frac {R{}}{n}\delta^k_j)
    \,
    d\sigma_k}_{=: (*)}
     \,.
           \label{2VIII21.2}
    \end{eqnarray}
We claim that the second line can be made arbitrarily small by choosing $i$ sufficiently large and $\epsilon$ sufficiently small.
For this we  apply the divergence theorem on the set
$$
  \red{\Omega_\epsilon}:= \{0\le x \le \epsilon/100 \}
  \times (\blue{\overline{D(2\epsilon)\setminus D(\epsilon/2)}}
  )
  \,.
$$
Letting
\begin{equation}\label{4VIII21.41}
  \myV^k:=  ( R{}^k{}_j - \frac {R{}}{n}\delta^k_j)  D^j ( \rphi_i V_1 )
    \,,
\end{equation}
where $\rphi_i$ has been extended away from $\partial M_1$ by requiring $\partial_r \rphi_i = 0$,
we have
\begin{equation}\label{4VIII21.42}
  \int_{\red{\Omega_\epsilon}} D_k \red{\myV}^k
   d\mu_g
   = (*)
   +
    \underbrace{
     \int_{  \red{\Omega_\epsilon}\cap \{x=\epsilon/100\}}\red{\myV}^k d\sigma_k
     }_{=:a}
     +
    \underbrace{
     \int_{[0,\epsilon/100]\times S_{2\epsilon}}
     \red{\myV}^k d\sigma_k
     }_{=:b}
     +
    \underbrace{
     \int_{[0,\epsilon/100]\times S_{ \epsilon/2}}
     \red{\myV}^k d\sigma_k
     }_{=:c}
   \,.
\end{equation}

On $[0,\epsilon/100]\times S_{ \epsilon/2}$ the metric is exactly hyperbolic, thus $c=0$.

For the purpose of Corollary~\ref{C21V22.2} below we will do the calculation that follows in general dimensions $n\ge 3$. We will assume \eqref{21V22.3} except that we require now
\begin{equation}\label{21V22.92}
  \sigma >  n-\frac1 2
  \,.
\end{equation}
Following \cite{ChDelayHPETv1} we define $s$ via the inequality
\begin{equation}\label{21V22.91}
 \mbox{$\sigma > (n-1)/2 +s$ for some
 $s\in (n/2,(n+1)/2)$}
\end{equation}
(actually $s=n/2$ is also allowed in \cite{ChDelayHPETv1} but we need $s>n/2$ in our calculations below).
We then have the following estimates in $\red{\Omega_\epsilon}$ (\cite{ChDelayExotic}, see also \cite[Remark~2.1]{ChDelayHPETv1}):
%
\begin{equation}\label{4VIII21.1}
  | g-b|_b+| D b|_b+| D^2 b|_b = O(\red{r^{-\sigma}}) + o (\epsilon^{\red{\sigma-s}}r^{-s})
 \,;
\end{equation}
recall that $D$ is the covariant derivative of $g$ and $r=1/x$.
 Let $\bar V = \rphi_i r \equiv \rphi_i V_1$. Then $d \bar V=\rphi_i dr+r\partial_A\rphi_i dx^A$  so that $|d\bar V|_b=\sqrt{r^2\rphi_i^2 +|d\rphi_i |^2_{h_0}}$.
We thus have
\begin{equation}\label{5VIII21.11a}
  |d\bar V|_b   \le \rphi_i r + |d\rphi_i |_{h_0}
  \,,
\end{equation}
and
\begin{equation}\label{5VIII21.11b}
  |R{}^i{}_j - \frac {R{}}{n}\delta^i_j |_b =
   \left\{
     \begin{array}{ll}
       0, & \hbox{in the hyperbolic region;} \\
   O(\red{r^{-\sigma}}) + o ({\epsilon}^{\red{\sigma-s}}r^{-s}), & \hbox{in the gluing region;} \\
        O(\red{r^{-\sigma}}), & \hbox{elsewhere.}
     \end{array}
   \right.
\end{equation}

We are ready now to consider the integral $b$ in \eqref{4VIII21.42}.
There we are in the last  case of \eqref{5VIII21.11b} so  that, setting  $x=1/r$, it holds
\begin{equation*}\label{4VIII21.51b}
   b =  \int_{[0,\epsilon/100]\times S_{2\epsilon} } O(x^{\sigma})
   \big( O(x^{-1})
    \rphi_i +  |d\rphi_i |_{h_0}
    \big) x^{1 -n}
    \, dx \,
    d^{n-2}\mu
		\,,
\end{equation*}
where $d^{n-2}\mu$ is the measure induced by $h_0$ on $S_{2\epsilon}$.
Convergence of the integral in $x$ requires $\sigma>n-1$ (which is satisfied in view of \eqref{21V22.92}), and after integration one obtains
\begin{equation}\label{4VIII21.51b2}
	b=
 \int_{ S_{2\epsilon} } O(\epsilon^{\sigma-n +1})
   \big(
    \rphi_i + \epsilon |d\rphi_i |_{h_0}
    \big)
    d^{n-2}\mu
\,.
\end{equation}

Consider, next, the integral $a$ in \eqref{4VIII21.42}.
At $r= 100/\epsilon$ we find, after taking into account the measure induced by $g$,
\begin{equation}\label{4VIII21.51}
   a =   \int_{\blue{\overline{D(2\epsilon)\setminus D(\epsilon/2)}} }
   \big( O(\red{\epsilon^{\sigma-n}})
    \rphi_i + O(\red{\epsilon ^{\sigma-n+1}}) |d\rphi_i |_{h_0}
    \big)
    d\mu_{h_0}
   \,.
\end{equation}

We pass now to the analysis of the integral over $\Omega_\epsilon$. A calculation shows that
$$
 \zD \zD  \bar V  =\bar V   b   + r \partial_A\partial_B \rphi_i \, dx^A   dx^B
$$
(recall that $\zD$ is the covariant derivative of the metric $b= dr^2/r^2+r^2 h_0$),
so that
$$
 |\zD \zD  \bar V  - \bar V   b|_b  = \frac{\sqrt{h_0^{AB} h_0 ^{CD} \partial_A\partial_C \rphi_i \partial_B\partial_D \rphi_i }} r
 \equiv \frac{|\partial^2 \rphi_i |_{h_0}}{r}
 \,.
$$
 Thus, with the error terms for tensors in $b$-norm unless explicitly indicated otherwise,
\begin{eqnarray}
 D_k D_j \bar V
  & = &  \zD_k\zD_j \bar V +  \big ( O(\red{r^{-\sigma}}) + o (\epsilon^{\red{\sigma-s}}r^{-s}) \big ) D \bar V
  \nonumber
\\
 &  =  &
   \bar V b_{kj}
    + O ( {|\partial^2 \rphi_i |_{h_0}} r^{-1})
   +  \big ( O(\red{r^{-\sigma}}) + o (\epsilon^{\red{\sigma-s}}r^{-s}) \big ) D \bar V
  \nonumber
\\
 &  =  &
   \bar V \big (g_{kj}+O(\red{r^{-\sigma}}) + o (\epsilon^{\red{\sigma-s}}r^{-s})\big  )
    + O ( {|\partial^2 \rphi_i |_{h_0}} r^{-1})
    \nonumber
\\
 &&   +  \big ( O(\red{r^{-\sigma}}) + o (\epsilon^{\red{\sigma-s}}r^{-s}) \big )
 \big(\rphi_i r + |d\rphi_i |_{h_0}
 \big)
  \nonumber
\\
 &  =  &
   \bar V   g_{kj}
    + O (r^{-1}) {|\partial^2 \rphi_i |_{h_0}}
     + \rphi_i \big( O(r^{\red{1-\sigma}}) + o (\epsilon^{\red{\sigma-s}}r^{1-s}) \big )
    \nonumber
\\
 &&
  + |d\rphi_i |_{h_0}  \big( O(\red{r^{-\sigma}})   + o (\epsilon^{\red{\sigma }-s}r^{ -s}) \big )
  \, .
  \label{5VIII201.1}
\end{eqnarray}
Next, the curvature scalar $R$ is constant, which implies that $D_i R^i{}_j=0$ by the
(twice contracted) second Bianchi identity. We thus find
\begin{eqnarray}
  D_k \red{\myV}^k
   & = &
    \underbrace{( R{}^k{}_j - \frac {R{}}{n}\delta^k_j)
    }_{
    O(\red{r^{-\sigma}}) + o (\epsilon^{\red{\sigma-s}}r^{-s})
    }
     D_k   D^j  \bar  V
 \nonumber
\\
 & = &
 (O(\red{r^{-\sigma}}) + o (\epsilon^{\red{\sigma-s}}r^{-s}))(O(r^{\red{1-\sigma}}) + o (\epsilon^{\red{\sigma-s}}r^{-s+1}))\rphi_i \nonumber
\\
 & & +
 (O(\red{r^{-\sigma}}) + o (\epsilon^{\red{\sigma-s}}r^{-s}))(O(\red{r^{-\sigma}}) + o (\epsilon^{\red{\sigma-s}}r^{-s }))
  |\red{d\rphi_i}| _{h_0}
  \nonumber
\\
 &    &
 +
    \big( O (r^{\red{-\sigma-1}})
   + o (\epsilon ^{\red{\sigma}- s} r^{- s-1})
   \big)
   |\partial^2 \rphi_i |_{h_0}
 \nonumber
\\
 & = &
 (O(r^{\red{-2\sigma+1}}) + o (\epsilon^{\red{\sigma-s}}r^{\red{-s-\sigma+1}}) + o (\epsilon^{\red{2\sigma-2s}}r^{-2s+1}))
  \rphi_i \nonumber
\\
 & & +
(O(r^{\red{-2\sigma}})
 + o (\epsilon^{\red{\sigma-s}}r^{\red{-s-\sigma}}) + o (\epsilon^{\red{2\sigma-2s}}r^{-2s }))
  |\red{d\rphi_i}| _{h_0}
  \nonumber
\\
 &    &
 +
    \big( O (r^{\red{-\sigma-1}})
   + o (\epsilon ^{\red{\sigma- s}} r^{- s-1})
   \big)
   |\partial^2 \rphi_i |_{h_0}
    \,.
    \label{4VIII21.410}
\end{eqnarray}
The requirement of convergence of the  integral in $r$ together with \eqref{21V22.91} leads to the restriction  $n<5$. Assuming this, we find
\begin{equation}\label{4VIII21.52}
     \int_{\red{\Omega_\epsilon}} D_k \red{\myV}^k
   d\mu_g
   = \int_{\blue{\overline{D(2\epsilon)\setminus D(\epsilon/2)}} }
    \big(
     O( \epsilon^{\red{2\sigma-n}} )
     \rphi_i +
     O( \epsilon^{\red{2\sigma-n+1}})
     |\red{d\rphi_i}| _{h_0}
     +
     O( \epsilon^{\red{\sigma-n+2}} )
      |\partial^2 \rphi_i |_{h_0}
      \big)
       d\mu_{h_0}
    \,.
\end{equation}
Collecting terms in \eqref{4VIII21.42} leads to the following form of \eqref{2VIII21.2}
\begin{eqnarray}
 \red{m_{i,\epsilon}}
  &  = &
  \hat m_{i,\epsilon}
+
 \int_{ S_{2\epsilon} } O(\epsilon^{\sigma-n +1})
   \big(
    \rphi_i + \epsilon |d\rphi_i |_{h_0}
    \big)
    d^{n-2}\mu
    \nonumber
\\
 &&
    +  \int_{\blue{\overline{D(2\epsilon)\setminus D(\epsilon/2)}} }
    \big(
     O(  \red{\epsilon^{\sigma-n}} )
     \rphi_i +
     O(  \red{\epsilon^{\sigma-n+1} } )
     |\red{d\rphi_i}| _{h_0}
     +
     O(  \red{\epsilon^{\sigma-n+2}} )
      |\partial^2 \rphi_i |_{h_0}
      \big)
       d\mu_{h_0}
    \,.
     \phantom{xxx}
           \label{15VIII21.1}
    \end{eqnarray}

Let $\eta>0$ and
  $\rphi_\infty = e^{-\mv_\infty/2}$. If $n=3$  (thus $\sigma> 5/2$), by \eqref{2VIII31.1}-\eqref{2VIII31.3}, \eqref{22V221}
    and \eqref{25VIII21.1}  we can choose $\epsilon$ small enough so that
\begin{eqnarray}
\lefteqn{
 \int_{ S_{2\epsilon} } O(\epsilon^{\sigma-n +1})
   \big(
    \rphi_i + \epsilon |d\rphi_i |_{h_0}
    \big)
    d^{n-2}\mu
}
&&
    \nonumber
\\
 &&
 +
 \int_{\blue{\overline{D(2\epsilon)\setminus D(\epsilon/2)}} }
    \big(
     O( \red{\epsilon^{\sigma-n}} )
     \red{\rphi_\infty} +
     O(  \red{\epsilon^{\sigma-n+1}}  )
     |\red{d\rphi_\infty}| _{h_0}
     +
     O(  \red{\epsilon^{\sigma-n+2}} )
      |\partial^2 \red{\rphi_\infty} |_{h_0}
      \big)
       d\mu_{h_0}
        < \eta/2
    \,.
 \phantom{xxxxxx}
           \label{15VIII21.1asdf}
    \end{eqnarray}
The function $\rphi_i$ tends uniformly on
$\blue{\overline{D(2\epsilon)\setminus D(\epsilon/2)}}$ to $\rphi_\infty $ when $i$ tends to infinity. This, together with \eqref{15VIII31.2}, shows that we can  choose   $i$ large enough so that
\begin{eqnarray}
\lefteqn{
 \int_{ S_{2\epsilon} } O(\epsilon^{\sigma-n +1})
   \big(
     \red{|\rphi_\infty-\rphi_i|} + \epsilon  \red{|d\rphi_\infty-d\rphi_i|}_{h_0}
    \big)
    d^{n-2}\mu
 }
 &&
    \nonumber
\\
&
 +
\displaystyle
 \int_{\blue{\overline{D(2\epsilon)\setminus D(\epsilon/2)}} }
    \Big(
&
     O(  \red{\epsilon^{\sigma-n}} )
     \red{|\rphi_\infty-\rphi_i|}  +
     O(  \red{\epsilon^{\sigma-n+1}}  )
     |\red{d(\rphi_\infty-\rphi_i)}| _{h_0}
    \nonumber
\\
 &&
     +
     O(  \red{\epsilon^{\sigma-n+2}} )
      |\partial^2 \red{(\rphi_\infty-\rphi_i)} |_{h_0}
      \Big)
       d\mu_{h_0}
        < \eta/2
    \,.
           \label{15VIII21.1asdmf}
    \end{eqnarray}
We conclude that the integral over the annulus  $\blue{\overline{D(2\epsilon)\setminus D(\epsilon/2)}} $ tends to zero as $\epsilon$ tends to zero and $i$ tends to infinity.

The above  analysis of the mass can be summarised as follows:

\begin{Proposition}
  \label{P21V22.1}
  Let $M$ have dimension three or four, and suppose that for every $\eta>0$ we can find $i$ large enough so that \eqref{15VIII21.1asdf}-\eqref{15VIII21.1asdmf} hold. Then, when $\epsilon$ tends to zero and $i$ tends to infinity the mass of each summand of $(M,g)$ tends to
\begin{eqnarray}
   -
   \lim_{x\rightarrow0}\int_{\{x\} \times \partial M_1}
   D^j ( \rphi_\infty  V_1 )
    ( R{}^k{}_j - \frac {R{}}{n}\delta^k_j)
    \,
    d\sigma_k
     \,,
           \label{21V22.11}
\end{eqnarray}
where the limit exists and is finite.
\end{Proposition}

By Remark~\ref{R22V22.1} and Deligne-Mumford compactness (cf., e.g., \cite[Proposition~A.2, Appendix~A.1]{RupflinToppingZhu}) the metric $e^{\omega_\infty}h_0$  is the hyperbolic metric on the punctured torus.
 This finishes the proof of Theorem~\ref{T29VII21.1}.
\qed

\begin{Remark}
  \label{R25V22.1}
{\rm
For further reference we comment on the possibility that the discrete family of circles $S_{1/i}$ across which the doubling has been done is replaced by  a
continuous family $S_\apar$, with $\apar\in (0,\apar_0)$ for some $\apar_0>0$, with associated family of solutions $\omega_{\apar}$. The only place in our arguments above where the discrete character of the parameter $\apar$ matters is the extraction of a diagonal subsequence so that the solutions converge to a limiting function. For this
we take any sequence $a_i$ converging to zero, so that the resulting sequence $\omega_{\apar_i}$ has a subsequence $\omega_{ {\apar_{i_j} }}$ converging to a limiting function, say
$\omega_{\{\apar_{i_j}\}}$, on compact subsets of $\T^2\setminus \{p\}$. By what has been said the metric  $e^{\omega_{\{\apar_{i_j}\}}} h_0$ is the unique hyperbolic metric on $\T^2\setminus\{p\}$ with puncture at $p$, hence is independent both of  the sequence ${\apar_{i_j}}$  and of the sequence ${\apar_{i }}$. Setting
 $ \omega_{\infty}:=  \omega_{\{\apar_{i_j}\}}$, it is then standard to show that $\omega_\apar$ converges to $\omega_\infty$ on any compact subset of $\T^2\setminus\{p\}$ as $\apar$ tends to zero, and  depends continuously on $\apar$ on any compact subset of $\T^2\setminus\{p\}$ in $C^\infty$, with  the length of the minimising geodesic in the middle of the connecting neck depending continuously upon $\apar$, and with the mass of $(M,g)$ depending continuously upon $\apar$.
}
\qed
\end{Remark}

It is of interest to enquire about the shape of the glued boundary-manifold. From what has been said the metric on the conformal boundary converges (uniformly on compact sets away from the puncture) to the cusp metric on  each of the two copies of  the punctured torus (cf. Figure~\ref{F21XI21.1}), with asymptotic behaviour near the puncture approximated by
(compare~\cite{MelroseZhu})
\begin{equation}\label{26VIII21.31}
  \frac{1}{\rho^{ 2}\red{\log}^2(\rho)} (dx^2 + dy^2) =  \frac{1}{\rho^{ 2}\red{\log}^2(\rho)} (d\rho^2 + \rho^2 d\varphi^2)
  \,.
\end{equation}
Since
$$
\int \frac 1 {\rho \red{\log} \rho} d\rho = \red{\log}( -\red{\log} \rho)\to_{\rho \to 0}\infty
\,,
$$
the connecting necks become  longer and longer, with a  circumference $\approx 2 \pi /|\red{\log} \rho| = 2 \pi /\red{\log} i $ tending to zero
when $i\to\infty$; see Figure \ref{F21XI21.2}.
\begin{figure}[b]
	\centering
 \includegraphics{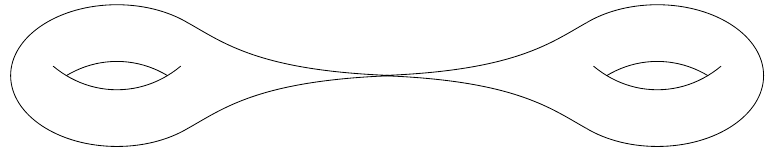}
    \caption{Two punctured tori, connected by a thin neck. In the limit $i \to \infty$
    the necks become longer and longer with circumferences shrinking to zero.}
	\label{F21XI21.2}
\end{figure}

\subsection{Higher genus}
 \label{ss1821.2}

It should be clear that one can iterate the construction of the last section to obtain three-dimensional conformally compact ALH manifolds without boundary, with constant scalar curvature, negative mass, and a conformal infinity of arbitrary genus. One possible way of doing this proceeds as follows:  Let $(M_1,g_1)$ be a three-dimensional ALH metric with constant scalar curvature and toroidal boundary at infinity. Let $N\in
\N$ and let $\{p_i\}_{i=1}^N$ be a collection of points lying on the conformal boundary $\partial M_1$ of $M_1$. Let $\mcU_i$ be any pairwise-disjoint family of coordinate half-balls
of coordinate radius $\epsilon_i$
centered at the points $p_i$.  Let $(M_2,g_2)$ be another ALH manifold with the same asymptotic behaviour and identical conformal geometry at infinity; e.g., an identical copy of $(M_1,g_1)$. The manifold $(M,g)$ is taken to be a boundary-gluing of $(M_1,g_1)$ with $(M_2,g_2)$, where each $\mcU_i\subset M_1$ is glued with its partner in $M_2$ in a symmetric way; Figure~\ref{26VIII21.1} illustrates what happens at the conformal boundary at infinity.
 \begin{figure}
	\centering
 \includegraphics[width=.5\textwidth]{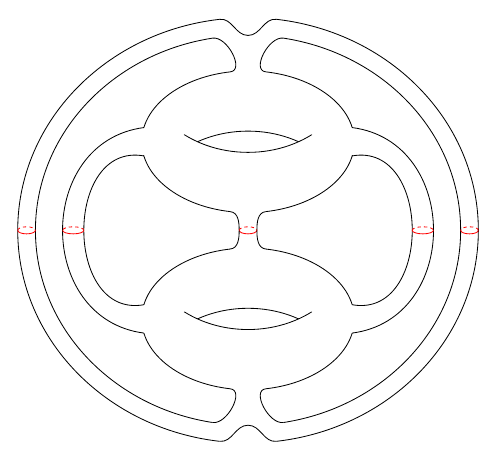}
    \caption{The gluing is done at the discs, whose boundaries become closed geodesics on the boundary, indicated in red in the figure. The construction is symmetric under the reflection across the horizontal plane passing through the center of the figure.  }
	\label{26VIII21.1}
\end{figure}
 Solving the Yamabe equation on $\partial M$,
for $j\in \N$ one obtains a function $\mv_j$ on
$$
\hat M_j  := \partial M_1 \setminus
\big( \cup_{i=1}^N  D(p_i, 1/j)
\big)
$$
such that the constant-Gauss-curvature representative of the conformal metric on $\partial M$ equals
$$
e^{\mv_j} h_0
$$
on $\zMone  $, with $\mv_j$ satisfying the Neumann boundary condition \eqref{1VIII31.3} on each of the coordinate circles $S(p_i)_{1/j}$ centred at $p_i$:
\begin{equation}\label{1VIII31.3a}
 \partial_\rho   \mv_j |_{S_{1/j}} = -2j
  \,.
\end{equation}
This leads to the integral identity
\begin{equation}\label{1VIII31.4a}
  \int_{\hat M_j} e^{\mv_j} d\mu_{h_0} = 2 N \pi
  \,,
\end{equation}
The analysis of Section~\ref{ss18VIII21.1} applies near each of the punctures $p_i$. Letting  $m(\epsilon_1,\cdots,\epsilon_{N},j)$
 denote the mass of $(M,g)$ we conclude that if we started with a Horowitz-Myers metric we will have
\begin{equation}\label{18VIII21.1a}
 \red{m(\epsilon_1,\cdots,\epsilon_{ N},j)}
 < 0
\end{equation}
%
for all $j$ larger than $j_0$, for some $j_0$ that depends upon $\max \epsilon_i$.

From \eqref{1VIII31.4a} we find
\begin{equation}\label{1VIII31.4ax}
 |\partial M|_{h_{-1}} = 2 \int_{\hat M_j} e^{\mv_j} d\mu_{h_0} = 4 N \pi
  \,.
\end{equation}
Equivalently, since the scalar curvature $R_{h_{-1}}$ of the metric $h_{-1}$ equals $-2$ we have
\begin{equation}\label{1VIII31.4axm}
 4\pi \chi(\partial M)= \int_{\partial M}  R_{h_{-1}} d\mu_{h_{-1}}= -4\int_{\hat M_j} e^{\mv_j} d\mu_{h_0} = -8 N \pi
  \,,
\end{equation}
which shows that $\chi(\partial M) = - 2 N$, and thus $\partial M$ has genus $N+1$.

A variation of the above is the following: Let $\psi: \partial M_1 \to \partial M_1$ be an isometry of $(M_1,h_0)$ such that $\psi\circ \psi$ is the identity map. Let $N$ be even, say $N=2
\hat N$, and suppose that the collection $\{p_i\}_{i=1}^N$ of $N$ distinct points  on $\partial M_1$ takes the form $\{q_i, \psi(q_i)\}_{i=1}^{\hat N}$. We can then pairwise glue, as above, neighborhoods of the points $q_i$ with  neighborhoods of their $\psi$-symmetric partners $\psi(q_i)$. Thus $\partial M$ is obtained by adding to $\partial M_1$ a family of $\hat N$ necks connecting the points $q_i$ and $\psi(q_i)$; see Figure~\ref{F26VIII21.3}.
 \begin{figure}
	\centering
 \includegraphics[width=.5\textwidth]{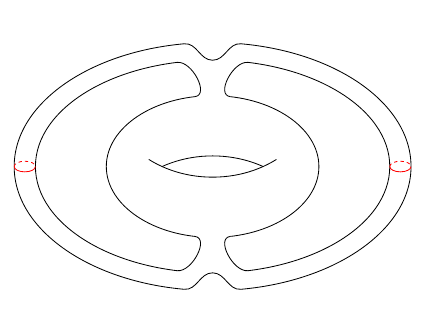}
    \caption{The symmetry $\psi$ is chosen to be a reflection across the vertical plane passing through the center of the figure.  }
	\label{F26VIII21.3}
\end{figure}

 Let $M_1$ be the space-part of the manifold \red{underlying} the Horowitz-Myers metric. If we choose $\psi$ and the necks in a way compatible with the constructions described in Section~\ref{ss18VIII21.1}, one obtains a family of ALH manifolds with constant scalar curvature,
 with conformal infinity of genus $\hat N +1$,
with masses $\red{m(\epsilon_1,\cdots,\epsilon_{\hat N},j)}$
which are negative for $j$ large enough.

\section{A lower bound on mass?}
 \label{s19VIII21}

It has been conjectured~\cite{HorowitzMyers,WoolgarRigidityHM} that the Horowitz-Myers metrics have the lowest mass within the class of ALH manifolds without boundary and with the same conformal metric at their toroidal infinity.
If this were wrong, then  the mass of conformally compact three-dimensional manifolds with constant scalar curvature and toroidal infinity could
 be arbitrarily negative at fixed conformal class at infinity. Hence our construction would provide higher-genus constant-scalar-curvature metrics which have arbitrarily negative mass at fixed conformal class at infinity.

As such, one can envision a construction which cuts a metric with conformal infinity of higher genus into pieces,  producing a finite-number of ALH metrics with toroidal infinity which would assemble together to the original one by a boundary-gluing in the spirit of this work. The validity of the Horowitz-Myers conjecture would then provide a lower bound for the mass of each summand.
It is thus tempting to formulate the following conjecture

\begin{conjecture}
 \label{C25VIII21}
Let $(M,g)$ be an $n$-dimensional conformally compact ALH manifold,  without boundary  {with well defined mass $m$} and with scalar curvature satisfying $R\ge -n(n-1)$. There exists a constant $m_0$,  which  depends only upon the conformal class of the metric at infinity  {and which is negative if this conformal class is not that of the round sphere}, such that
\begin{equation}\label{21VIII21.11}
  m \ge m_0
  \,.
\end{equation}
\end{conjecture}

It is also tempting to conjecture existence of  {unique} static metrics  as above  which saturate the inequality if $m_0$ is optimal.

\appendix

\section{A punctured torus}
 \label{sApp20XI21.1}
 Recall that the model for the cusp (at $t=-\infty$ or $\rho=0$) is :
 $$dt^2+e^{2t}d\varphi^2=\frac{1}{\rho^2\red{\log}^2(\rho)}(d\rho^2+\rho^2d\varphi^2),$$
 with $e^{-t}=-\red{\log}(\rho)$. Here the Gauss curvature is $-1$, thus the Ricci scalar equals $-2$.

 The simplest hyperbolic one-punctured  torus, also called ``the hyperbolic torus with one cusp'',  is $(\R^2\backslash \Z^2)/\Z^2$ with a hyperbolic metric.
 Let us use the notation of \cite{Sugawa}
 (keeping in mind a different scaling of the hyperbolic metric there)
 in order to compare the hyperbolic metric and the flat metric on  this torus.
 For a domain $D$ in $\C$, the complete hyperbolic  metric on $D$ with Gauss curvature equal to $-1$,  as obtained using the Riemann mapping theorem applied to the universal covering space, is denoted by
 $4 \rho_D^2(z)|dz|^2$. If $D\subset D'$ then  by the comparison principle \eqref{14IX21.2as}, also known as the Nevanlinna principle, we have $\rho_D\geq\rho_{D'}$.
 Let us set
 $$
 D_0=\hat\C\backslash \{0,1,\infty\}\;,\;\;\mathbb D^*=B(0,1)\backslash \{0\}\;,\;\;\Z[i]=\{m+in\;, m,n\in\Z\}=\Z^2\;,\;\; \tilde D=\C\backslash \Z[i]
 \,.
 $$
 Because $\mathbb D^*\subset\tilde D\subset D_0$, by the Hempel inequality~\cite{Hempel} (see, e.g., \cite{Sugawa}, inequality (1.2)) we have
 \begin{equation}\label{20XI21.11}
  \frac{1}{|z|(|\log(|z|)|+C_1)}\leq
  2 \rho_{D_0}(z)\leq
   2 \rho_{\tilde D}(z)\leq
  2
   \rho_{\mathbb D^*}(z)=\frac{1}{|z||\log(|z|)|}
 \,,
 \end{equation}
 where $C_1=\Gamma(1/4)^4/4\pi^2\approx 4.37688$.

 Let $p:\C\longrightarrow T_0=\C / \Z[i]$ be the canonical projection and denote by $[z]=p(z)$ the equivalence class of $z$.
 Let $D=T_0\backslash\{[0]\}$, then $D$ is proper subdomain of $T_0$ and $p^{-1}(D)=\tilde D$. On $D$,
 $\rho_D([z])=\rho_{\tilde D}(z)$ is well defined.  This defines the ``hyperbolic density'' $\rho_{X_0}$  on $X_0=\tilde D/ \Z[i]$, which is the simplest one-punctured hyperbolic torus.

 Because of the symmetries of $\tilde D$, the inequalities (\ref{20XI21.11}) on $\{z := x+iy, \;\, x,y\in[-1/2,1/2]\}\backslash \{0\}$ provide an estimate of  $\rho_{X_0}$ throughout $X_0$.
 It follows from \eqref{20XI21.11} that near the origin we have
 $$
 2 \rho_{\tilde D}=\frac{1}{|z||\log(|z|)|}\left(1+O\left(\frac1{|\log(|z|)|}\right)\right)
  \,.
 $$

 Note also that by \cite{Sugawa} the minimum of $\rho_{X_0}$ is attained at  $\frac{1+i}2$ and is given by his formula (1.3):
 \begin{equation}\label{20XI21.9}
   \inf 2 \rho_{X_0} = \frac{4 \pi ^{3/2}}{\Gamma(1/4)^2} \approx   1.694
   \,.
 \end{equation}

\section{Static KIDs}
\label{A17VIII21.1}
In this appendix we find all static potentials (KIDs) for the Birmingham-Kottler metric and the Horowitz-Myers
metric in
$(n+1)$-dimensions.
Let us summarise the argument.
We start by computing the KIDs for the Birmingham-Kottler
metric \eqref{genKott}.
We find that for $\mc \neq 0$ and $n \neq 2$ the solution to the KID equations is given by \eqref{VBKgen}.
Otherwise, the solution to the
KID  equations takes the form \eqref{v1}. For $\mc = 0$ and $n \neq 2$
 the function
 $\beta$ of \eqref{v1}
 is of the form \eqref{betamc0k} or \eqref{betamc0k0} and $\Omega$
is a solution to the differential equations \eqref{Omegaequation}.
In the case $n = 2$, the function
$\beta$ is of the form \eqref{betan21} or \eqref{betan22} and $\Omega$ takes the form \eqref{24XI21.1} .
Next, we consider the Horowitz-Myers metric \eqref{HMmetric}, which for $\mc = 0$ is of Birmingham-Kottler type. For $\mc \neq 0$ we show
that the solution to the KID  equation
is given by \eqref{HMmnot0final}.

\subsection{Birmingham-Kottler}
\label{BKkids}
The Birmingham-Kottler metric in $(n+1)$-dimensions, $n\ge 2$, which we denote by $\tilde{g}$, reads
\begin{equation}
\label{genKott}
   \tilde{g} = - f(r) dt^2 + \frac{1}{f(r)} dr^2 + r^2 \mathring h_{A B} dx^A dx^B\,,
\end{equation}
where $\mathring h_{A B} = \mathring h_{A B}(x^C)$ where $A,B =1, ... n-1$, and with
\begin{equation}
    f(r) =  r^2 + k - \frac{2 \mc}{r^{n-2}}
\end{equation}
and
\begin{equation}
    R(\mathring h) = k (n-1) (n-2)
    \,,\quad
     k \in \{0,\pm 1
      \}
     \,.
\end{equation}
For $n = 2$, $R(\mathring h)$ vanishes and thus we set $k = 0$ without loss of generality. This then yields that
\begin{equation}
    f(r) = r^2 - 2 \mc\,
\end{equation}
for $n =2$.
The metric \eqref{genKott} fulfills the Einstein  equations with negative
cosmological constant
\begin{equation}
    \Lambda = - \frac{n (n-1)}{2 }\,.
\end{equation}
Using this normalisation of $\Lambda$, the Einstein equations read
\begin{equation}
   {R}_{a b}(\tilde g) = - n {\tilde g}_{a b}\,, \qquad {R(\tilde g)} =  - n (n+1)\,,
\end{equation}
where $a, b$ run from $1, ..., n+1$.
To compute the KIDs we need the Ricci scalar and the Ricci tensor of the spatial part of $\tilde g$,  which we denote
 as $g$ in
what follows. These can be obtained from the Gauss-Codazzi equation \cite{Wald:book}
\begin{equation}
    R_{a b c d} (g) =  R^{ {\sparallel}}_{a b c d} (\tilde g) - K_{a c} K_{b d} + K_{b c} K_{a d}\,.
\end{equation}
where $K_{a b} $ is the extrinsic curvature of the level sets of $t$ (zero in our case)
\begin{equation}
  K_{a b}  = g_a {}^c  {\tilde D}_{c} \tilde n_{b}\,, 
\end{equation}
with $\tilde n_a$ being the normal one-form with the normalization chosen such that $\tilde n^a \tilde n_a = -1$.
The superscript ``$\red{\sparallel }$'' means that the
respective tensor is projected to the submanifold, e.g.
$R^{\red{\sparallel }}_{a b c d} (\tilde g) = g_a{}^e g_b{}^f g_c{}^g g_d{}^h R_{e f g h} (\tilde g)$.
Contracting with $g^{ac}$ we obtain
\begin{equation}
    R_{a b} (g) = R^{\red{\sparallel }}_{a b}(\tilde g) +\left(R_{c a d b} (\tilde g) \tilde n^c \tilde n^d \right)^{\red{\sparallel }} - K K_{a b} + K_{a c} K^{c}{}_b
\end{equation}
with \cite[p.\ 518]{MTW}
\begin{equation}
    \left(R_{c a d b} (\tilde g) \tilde n^c \tilde n^d \right)^{\red{\sparallel }} =
    - \mathcal L_{\tilde n} K_{a b} + K_{a c} K^c{}_b
    +  D_{(a} a_{b)} + a_a a_b
\end{equation}
where $a_a = \tilde n^b  {\tilde D}_{b} \tilde n_a$, $D_{a}$ denotes the covariant derivative of the metric $g$ and symmetrization is defined with
a factor $1/2$ such that
\begin{equation}
    \label{Riccisubmanifold}
    R_{a b} (g) = R^{\red{\sparallel }}_{a b}(\tilde g)
    - \mathcal L_{\tilde n} K_{a b} + 2 K_{a c} K^{c}{}_b- K K_{a b}
    +  D_{(a} a_{b)} + a_a a_b
    \,.
\end{equation}
The normal one-form to $t = \mathrm{const.}$ surfaces reads
\begin{equation}
\tilde n = \sqrt{f} dt
\end{equation}
such that
$
{\tilde n}_a {\tilde g}^{a b} {\tilde n}_b = - 1\,.
$
Finally
\begin{equation}
\tilde g_{a b} = g_{a b} - \tilde n_a \tilde n_b
\,,
\end{equation}
yielding that in coordinates $t, r, x^A$, the tensor field $g$ reads
\begin{equation}
g_{a b} =
\begin{pmatrix}
0 & 0 & 0\\
0 & 1/f(r) & 0 \\
0 & 0 & r^2 \mathring h_{A B}
\end{pmatrix}\,, \qquad g^{a}{}_b =
\begin{pmatrix}
0 & 0 & 0\\
0 & 1 & 0 \\
0 & 0 & 1_{A B}
\end{pmatrix}
 \,.
\end{equation}
 The only-nonvanishing Christoffel symbols (up to symmetries),
\begin{equation}
\tilde\Gamma^{c}{}_{a b} = \frac{1}{2}\tilde g^{c d} \left(
\tilde \partial_a\tilde g_{bd} + \tilde \partial_b\tilde g_{a d} - \tilde \partial_d\tilde g_{a b}
\right)
\,,
\end{equation}
are
\begin{align}
\tilde \Gamma^t{}_{t r} &= \frac{\partial_r f}{2 f} \,, \quad
\tilde \Gamma^r{}_{r r} = - \frac{\partial_r f}{2 f} \,, \quad
\tilde \Gamma^r{}_{t t} = - \frac{f \partial_r f}{2} \,, \quad
\tilde \Gamma^r{}_{A B} = - r \mathring h_{A B} f(r)\,,  \\
\tilde \Gamma^C{}_{r A}  
 &= \frac{\delta^{C}{}_A}{r} \,, \quad
\tilde \Gamma^C{}_{A B} = \frac{1}{2} \mathring h^{C D}
\left( \partial_A \mathring h_{BD} + \partial_B \mathring h_{A D} - \partial_D \mathring h_{A B}
\right)\,.
\end{align}
As already mentioned, in our case all components of $K_{ab}$ vanish,
so that the expression of the Ricci tensor of the constant $t$ submanifold \eqref{Riccisubmanifold}
 simplifies to
\begin{equation}
R_{a b}(g) = {R}_{a b}^{\red{\sparallel }}(\tilde g) +  D_{(a} a_{b)} + a_a a_b\,.
\end{equation}
We have that $a = a_r dr$ with
\begin{equation}
a_r = \tilde n^t {\tilde D}_{t} \tilde n_r = \frac{1}{2 f} \partial_r f
\end{equation}
and
\begin{align}
D_{r} a_r &= - \frac{1}{4 f^2} (\partial_r f)^2 + \frac{1}{2 f} \partial_r \partial_r f\,, &
D_A a_B &= \frac{\partial_r f}{2} r \mathring h_{AB}\,.
\end{align}
The only non-vanishing components of the Ricci tensor of the metric induced on the submanifolds of  constant $t$ are
\begin{align}
R_{r r} & = - \frac{n}{f}
+ \frac{1}{2 f} \partial_r \partial_r f\,,  &
R_{A B} &= \left(  - n r^2 + \frac{\partial_r f}{2} r \right) \mathring h_{AB}
\end{align}
and thus
\begin{equation}
R = R_{rr} g^{rr} + R_{A B} g^{A B} = - n (n-1)\,,
\end{equation}
where we have used the explicit form of $f(r)$.
Now we consider the KID  equation,
\begin{equation}
    A_{i j} := {D}_i {D}_j V - V \left({R}_{i j}- \frac{{R}}{n-1} {g}_{i j}\right)= {D}_i {D}_j V - V ({R}_{i j}+ n {g}_{i j}) = 0
    \,,
\end{equation}
where $i, j \in 1, ... n-1$. We have that
\begin{align}
A_{r r} &= \partial_r \partial_r V + \frac{1}{2 f} \partial_r f \partial_r V - \frac{V}{2 f}\partial_r \partial_r f = 0\,, \\
A_{r A} &= \partial_r \partial_A V - \frac{\delta_A{}^C}{r} \partial_C V
= \partial_r \partial_A V - \frac{1}{r} \partial_A V  =0\,, \\
A_{A B} &= \partial_A \partial_B V - \Gamma^{C}{}_{A B} \partial_C V + r \mathring h_{A B} f \partial_r V
- \frac{V \partial_r f}{2} r \mathring h_{AB} = 0 \label{aAB}\,.
\end{align}
Integrating the second equation in $r$ we have that
\begin{equation}
\label{v1}
\boxed{
V(r, x^1, ..., x^{n-1}) = \beta(r) +r \Omega(x^1, ...x^{n-1})\,,
}
\end{equation}
for some functions $\beta$ and $\Omega$.
Plugging \eqref{v1} into $A_{r r}$ we find
\begin{equation}
\label{arrplugged}
\Omega(x^1, ...x^{n-1}) \frac{\left(\partial_r f - r \partial_r \partial_r f \right)}{f}
+ \frac{\partial_r f \partial_r \beta - \beta \partial_r \partial_r f + 2 f \partial_r \partial_r \beta}{f} = 0\,.
\end{equation}
Thus, the differential equation can only be fulfilled if \emph{either} $\Omega(x^1, ...x^{n-1})  $ is constant, which can then without loss of generality be taken to be zero as any constant can be reabsorbed into a
redefinition of $\beta(r)$,  \emph{or} if
\begin{equation}
    \partial_r f - r \partial_r \partial_r f  \equiv  2 \mc (n-2) n r^{1-n} = 0
\end{equation}
which can only be fulfilled if $\mc = 0$ or $n = 2$.

 The remaining part of  \eqref{arrplugged} reads
\begin{equation}
    \partial_r f \partial_r \beta - \beta \partial_r \partial_r f + 2 f \partial_r \partial_r \beta =0\,.
\end{equation}
Changing the dependent variable $\beta$ to
\begin{equation}
\beta(r) = \sqrt{f(r)} \beta_2(r)
\end{equation}
we obtain the differential equation
\begin{equation}
3 \partial_r f \partial_r  \beta_2 +2 f \partial_r \partial_r \beta_2 = 0\,.
\end{equation}
From this we obtain that
\begin{equation}
\beta_2(r) = \tilde c - \hat{c}  \int \frac{\,dr}{f(r)^{\frac{3}{2}}}
\,,
\end{equation}
which yields that
\begin{equation}
\label{betageneral}
\beta(r) =  \tilde c \sqrt{f(r)}- \hat{c} \sqrt{f(r)} \int \frac{d  r}{f( r)^{\frac{3}{2}}}   \,.
\end{equation}

\subsubsection{\texorpdfstring{$\mc = 0$, $n\ne 2$}{mc = 0, n != 2}}
In the case $\mc = 0$ the function $f$ reduces to
\begin{equation}
    \label{mc0}
    f(r) =  r^2 + k
     \,,
\end{equation}
and \eqref{betageneral} readily integrates to
\begin{equation}
    \label{betamc0k}
    \boxed{
\beta(r) =- \frac{\hat{c} r}{k} + \tilde c \sqrt{r^2 + k}
    }
\end{equation}
if $k \neq 0$, and to
\begin{equation}
    \label{betamc0k0}
    \boxed{
\beta(r) =  \frac{\hat{c}}{2 r} + \tilde c r
    }
\end{equation}
if $k = 0$.
The differential equations \eqref{aAB} reduce to
\begin{align}
    \label{Omegaequation}
    \boxed{
D_A D_B \Omega\equiv \partial_A \partial_B \Omega -\Gamma^{C}_{A B} \partial_C \Omega =
 \mathring h_{A B} \left(\hat c - k \Omega\right)\,.
    }
\end{align}
This is an overdetermined system of equations, and all triples  $(M, \mathring h,\Omega)$ satisfying \eqref{Omegaequation}, where $\Omega\not\equiv 0$ and $(M,\mathring h)$ a complete Riemannian manifold, were
classified in \cite{Tashiro}.

In the special case
\begin{equation}
    \mathring h_{A B} dx^A dx^B = d\ntheta^2 + {(d\blue{\theta}^1)}^2 + ... + {(d\blue{\theta}^{n-2})}^2
\end{equation}
the constant $k$ is zero and
all Christoffel symbols of the boundary metric $\mathring h$ vanish, so that \eqref{Omegaequation} becomes
\begin{align}
\partial_A \partial_B \Omega= \hat{c} \,\delta_{A B}\,.
\end{align}
We then find
\begin{equation}
    \boxed{
    V(r) = \tilde c r+ \frac{\hat{c} }{2 r} + r \left(c \ntheta + \hat{c} \frac{\ntheta^2}{2}\right)
    + r\sum_{I = 1}^{n- 2}  c_I \theta^I + \hat{c} \frac{|\theta|^2}{2}  \,,
    }
\end{equation}
where $c$ and $c_I$ are 
 constants\,.

\subsubsection{\texorpdfstring{$\mc \neq 0$, $n\ne 2$}{mc != 0, n != 2}}
 \label{ss24XI21.1}

In the generic case $\mc \neq 0$ the function $\Omega(x^1, ...x^{n-1}) $ vanishes identitcally and the function $V$ reduces to
\begin{equation}
V(r, x^1, ... x^{n-1}) = \beta(r) = \tilde c\sqrt{f(r)} - \hat{c} \sqrt{f(r)} \int \frac{d r}{f(r)^{\frac{3}{2}}} \,.
\end{equation}
The final set of equations that needs to be solved is \eqref{aAB}, which in this case reduces to
\begin{equation}
A_{A B} = r \mathring h_{A B} f \partial_r V
- \frac{V \partial_r f}{2} r \mathring h_{AB} = 0\,.
\end{equation}
Contracting with $\mathring h^{A B}$ yields
\begin{equation}
r f \partial_r V
- \frac{V \partial_r f}{2} r  = 0
\end{equation}
which implies that $\hat{c} = 0$ and thus
\begin{equation}
    \label{VBKgen}
    \boxed{
V(r, x^1, ..., x^{n-1}) = \tilde{c} \sqrt{f(r)} \,.
    }
\end{equation}

\subsubsection{\texorpdfstring{$n = 2$}{n = 2}}
 \label{s13I22.1}

As already discussed, in the case $n =2$ the function $f(r)$ reduces to
\begin{equation}
    f(r) = r^2 - 2 \mc\,,
\end{equation}
which has the same functional form as \eqref{mc0} but with $k$ in \eqref{mc0} now being replaced by
$ - 2 \mc$.
Again, \eqref{betageneral} thus readily integrates to
 \begin{equation}
  \label{betan21}
\beta(r) = \frac{\hat{c} r}{ 2 \mc} + \tilde c \sqrt{r^2 - 2 \mc}
\end{equation}
if $\mc \neq 0$, and to
\begin{equation}
  \label{betan22}
\beta(r) =  \frac{\hat{c}}{2 r} + \tilde c r
\end{equation}
if $\mc = 0$.
 The differential
equations \eqref{aAB} reduce to just one differential equation as the boundary metric is one-dimensional,
\begin{align}
D_A D_A \Omega =
 \mathring h_{A A} \left(\hat c +2 \mc \Omega\right)
  \,,
\end{align}
with $A = 1$. In a coordinate system $x^1=\psi$ in which $\mathring h_{11}$ equals $1$ the solutions are
\begin{equation}\label{24XI21.1}
  \Omega = \left\{
             \begin{array}{ll}
               \frac{\hat c}2 \psi^2 + a \psi + b, & \hbox{$\mc =0$;} \\
               a e^{\sqrt{2\mc} \psi} +b e^{-\sqrt{2\mc} \psi}
 -   \frac{\hat c}{2 \mc} , & \hbox{$\mc >0$;} \\
               a \sin({\sqrt{-2\mc} \psi}) +b\cos({\sqrt{-2\mc} \psi} )
 -   \frac{\hat c}{2 \mc}  , & \hbox{$\mc <0$,}
             \end{array}
           \right.
\end{equation}
with  constants $a,b\in \R$.

\subsection{Horowitz-Myers}\label{s20VII21.1}

We consider the Horowitz-Myers metrics, which take  the form
\begin{equation}
    \label{HMmetric}
    ds^2 =  \frac{1}{f(r)} dr^2 +
     f(r)  d \ntheta^2 +r^2 \left(-  dt^2 + (d\blue{\theta}^1)^2 + (d\blue{\theta}^2)^2 + ... + (d\blue{\theta}^{n-2})^2 \right)
     \,,
\end{equation}
where
\begin{equation}
    f(r) = r^2 - \frac{2 \mc}{r^{n-2}}\,,
\end{equation}
with $n \ge 2$.
The spatial part of the metric,
which  reads
\begin{equation}
    \label{HMTKspatial}
g_{i j} dx^i dx^j= \frac{1}{f(r)} dr^2 +f(r)  d \ntheta^2 +r^2 \left((d\blue{\theta}^1)^2 + (d\blue{\theta}^2)^2 + ... + (d\blue{\theta}^{n-2})^2 \right)\,,
\end{equation}
will sometimes also be referred to as the HM metric.
The indices $i, j$ run from $1$ to $n$.

The only non-vanishing Christoffel symbols (up to symmetries) for \eqref{HMTKspatial}
are
\begin{subequations}
\begin{align}
    {\Gamma}^r{}_{rr} &= - \frac{\partial_r f}{2 f} \,, &
    {\Gamma}^r{}_{\ntheta \ntheta}  &= - \frac{f \partial_r f}{2}  \,,  &
    {\Gamma}^r{}_{\blue{\theta}^I \blue{\theta}^I}  &= - r f
    \,, &
    {\Gamma}^{\ntheta}{}_{\ntheta r} &= \frac{\partial_r f}{2 f}\,,  &
    {\Gamma}^{\blue{\theta}^I}{}_{\blue{\theta}^I  r} &= \frac{1}{r}\,. &
\end{align}
\end{subequations}
Next we compute
\begin{equation}
    {D}_i {D}_j V(r, \ntheta, \blue{\theta}^1,... \blue{\theta}^{n-2}) = \partial_i \partial_j V(r, \ntheta, \blue{\theta}^1,... \blue{\theta}^{n-2}) - {\Gamma}^k{}_{i j} \partial_k V(r, \ntheta, \blue{\theta}^1,... \blue{\theta}^{n-2})\,.
\end{equation}
 We have that
 \begin{subequations}
 \label{didjN}
 \begin{align}
     {D}_r {D}_r V &= \partial_r \partial_r V + \frac{1}{2 f} \partial_r f \partial_r V\,, &
     {D}_r {D}_\ntheta V &= \partial_r \partial_\ntheta V - \frac{1}{2} \frac{\partial_r f}{f} \partial_\ntheta V\,, \\
     {D}_r {D}_{\blue{\theta}^I} V &=  \partial_r \partial_{\blue{\theta}^I}  V -\frac{1}{r} \partial_{\blue{\theta}^I}  V\,, &
     {D}_\ntheta {D}_\ntheta V &= \partial_\ntheta \partial_\ntheta  V + \frac{1}{2} f \partial_r f \partial_r V\,, \\
     {D}_\ntheta {D}_{\blue{\theta}^I} V &= \partial_\ntheta \partial_{\blue{\theta}^I} V\,, &
     {D}_{\blue{\theta}^I} {D}_{\blue{\theta}^J} V &= \partial_{\blue{\theta}^I} \partial_{\blue{\theta}^J} V + r f \partial_r V \delta_{I J}\,,
 \end{align}
 \end{subequations}
 where $I, J$ run from from $1$ to $n-2$.
 The only non-vanishing components of the Riemann tensor,
\begin{equation}
    {R}_{i j k}{}^{l} = \partial_j {\Gamma}^l{}_{i k} - \partial_i {\Gamma}^l{}_{j k} + \left({\Gamma}^m{}_{i k} {\Gamma}^l{}_{m j} - {\Gamma}^m{}_{j k} {\Gamma}^l{}_{m i} \right)\,,
\end{equation}
read (up to the ones obtained from antisymmetry in the first index pair)
\begin{subequations}
\begin{align}
    {R}_{r \ntheta r}{}^{\ntheta} &=   - \frac{\partial_r \partial_r f}{2 f} \,, &
     {R}_{r \ntheta \ntheta}{}^{r} &=  - \frac{f \partial_r \partial_r f}{2}\,, &
     {R}_{r \blue{\theta}^I r}{}^{\blue{\theta}^I} &=  - \frac{\partial_r f}{2 f r}\,, &
     {R}_{r \blue{\theta}^I \blue{\theta}^I}{}^{r} &= \frac{ r\partial_r f}{2} \,, \\
     {R}_{\ntheta \blue{\theta}^I \ntheta}{}^{\blue{\theta}^I} &= - \frac{f \partial_r f}{2 r}\,, &
     {R}_{\ntheta \blue{\theta}^I \blue{\theta}^I}{}^{\ntheta} &=  \frac{r \partial_r f}{2}  \,, &
     {R}_{\blue{\theta}^I \blue{\theta}^J \blue{\theta}^I}{}^{\blue{\theta}^J} &=  - f&\mathrm{with}&\, I \neq J\,.
\end{align}
\end{subequations}
From this we obtain that the only non-vanishing components of the Ricci-tensor,
\begin{equation}
    {R}_{i j} =  {R}_{i k j}{}^{k}
    \,,
\end{equation}
 read
 \begin{subequations}
 \label{ricci}
\begin{align}
    {R}_{r r} &= - \frac{\partial_r \partial_r f}{2 f} - \frac{(n-2)}{2} \frac{\partial_r f}{r f} \,,&
    {R}_{\ntheta \ntheta} &= - \frac{1}{2} f \partial_r \partial_r f - \frac{(n-2)}{2} \frac{f \partial_r f}{r }\,,& \\
    {R}_{\blue{\theta}^I \blue{\theta}^I} &= - (n-3) f - r \partial_r f\,. & &
\end{align}
This leads to the following formula for the Ricci scalar
\begin{equation}
    R = - n (n-1)\,.
\end{equation}
The KID  equations,
\begin{equation}
    A_{i j} := {D}_i {D}_j V - V \left({R}_{i j}- \frac{{R}}{n-1} {g}_{i j}\right)= {D}_i {D}_j V - V ({R}_{i j}+ n {g}_{i j}) = 0
     \,,
\end{equation}
\end{subequations}
thus become
\begin{subequations}
   \begin{align}
     A_{rr} =& ~\partial_r \partial_r V + \frac{1}{2f} \partial_r f \partial_r V - V
     \left(- \frac{\partial_r \partial_r f}{2 f} - \frac{(n-2)}{2} \frac{\partial_r f}{r f}
      + \frac{n}{f}
     \right)= 0 \,, \label{eqrr}\\
     A_{r \ntheta} =&~A_{\ntheta r} =\partial_r \partial_\ntheta V - \frac{1}{2} \frac{\partial_r f}{f} \partial_\ntheta V = 0 \label{rtheta}\,,\\
     A_{r \blue{\theta}^I} =&~A_{\blue{\theta}^I r} =\partial_r \partial_{\blue{\theta}^I} V - \frac{1}{r} \partial_{\blue{\theta}^I} V = 0 \label{rthetaI}\,,\\
     A_{\ntheta \ntheta} =&~ \partial_\ntheta \partial_\ntheta V + \frac{1}{2} f \partial_r f \partial_r V - V \left( - \frac{1}{2} f \partial_r \partial_r f
     - \frac{(n-2)}{2} \frac{f \partial_r f}{r} + n f \right) =0 \label{eqthetatheta}\,,
   \\
     A_{\ntheta \blue{\theta}^{I}} = &~ \partial_\ntheta \partial_{\blue{\theta}^I} V=0\,, \label{thetathetaI}
   \\
     A_{\blue{\theta}^{I} \blue{\theta}^{J}} = &~ \partial_{\blue{\theta}^I} \partial_{\blue{\theta}^J} V
      + \delta_{I J} \Big(
     r f \partial_r V - V \bigl( - (n-3) f - r \partial_r f
       + n r^2\bigr)
       \Big)
       =0 \label{thetaIthetaImt}\,.
   \end{align}
   \end{subequations}
   Solving \eqref{thetathetaI} and \eqref{thetaIthetaImt} with $I\ne J$     we obtain that
\begin{equation}
V(r, \ntheta, \blue{\theta}^1, ... \blue{\theta}^{n-2}) = \alpha(r, \ntheta)+\alpha_1(r, \blue{\theta}^1)+ \alpha_2(r, \blue{\theta}^2)+... + \alpha_{n-2}(r, \blue{\theta}^{n-2})\,.
\end{equation}
Plugging this into the equation \eqref{rthetaI} we get that
\begin{equation}
    \alpha_I(r,\blue{\theta}^I ) = \beta_I(r) + r {\hat \gamma}_I(\blue{\theta}^I)\,.
\end{equation}
Solving \eqref{rtheta} yields
\begin{equation}
    \alpha(r, \ntheta) = \tilde \beta(r) + \sqrt{f(r)} {\hat \gamma}(\ntheta)\,.
\end{equation}
Setting
$\beta(r) = \tilde \beta(r) + \sum_{I =1}^{n-2} \beta_I(r)$
we have
\begin{equation}
    V(r, \ntheta, \blue{\theta}^1, ... \blue{\theta}^{n-2}) = \beta(r) + \sqrt{f(r)} {\hat \gamma}(\ntheta) + r \sum_{I =1}^{n-2} {\hat \gamma}_I(\blue{\theta}^I)
    \,,
\end{equation}
with the last term omitted when $n=2$.
From \eqref{thetaIthetaImt} we obtain (by pairwise substraction)
\begin{equation}
    \partial_{\blue{\theta}^I} \partial_{\blue{\theta}^I} V = \partial_{\blue{\theta}^J} \partial_{\blue{\theta}^J} V
\end{equation}
for any $I, J$, which is equivalent to
\begin{equation}
    \partial_I \partial_I \hat\gamma_I(\blue{\theta}^I)=  \partial_J \partial_J\gamma_J(\blue{\theta}^J)\,.
\end{equation}
From this we obtain that (the $\blue{\theta}^I$-independent integration constant can be absorbed into $\beta(r)$)
\begin{equation}
\label{TKHMNsimp}
   V(r, \ntheta, \blue{\theta}^1, ... \blue{\theta}^{n-2}) = \beta(r) + \sqrt{f(r)} {\hat \gamma}(\ntheta)
        + r  c_I \blue{\theta}^I
          + r \frac{\hat{c}}{2} \blue{|\theta|} ^2
          \,.
\end{equation}
We turn now our attention to the equation
\begin{align}
    A_{\ntheta \ntheta}-f(r)^2 A_{rr} &=0 = \partial_\ntheta \partial_\ntheta V- f^2 \partial_r \partial_r V \nonumber \\
    &= \frac{1}{4} \sqrt{f(r)} \left(\hat \gamma(\ntheta) \left((\partial_r f)^2 - 2 f \partial_r \partial_r f\right) - 4 f^{\frac{3}{2}} \partial_r \partial_r \beta + 4 \partial_\ntheta  \partial_\ntheta  \hat \gamma
    \right)\,,
\end{align}
from which we find
\begin{equation}
\label{ArrmAttHM}
    \hat \gamma(\ntheta) \left((\partial_r f)^2 - 2 f \partial_r \partial_r f\right) + 4 \partial_\ntheta  \partial_\ntheta  \hat \gamma = 4 f^{\frac{3}{2}} \partial_r \partial_r \beta \,.
\end{equation}
For the function $f(r) =r^2- \frac{2 \mc}{r^{n-2}}$ we have that
\begin{equation}
\label{HMbeforesplit}
    \left((\partial_r f)^2 - 2 f \partial_r \partial_r f\right) = 4 n \mc r^{2-2 n} \left((-1+n) r^n - (-2+n) \mc \right)\,.
\end{equation}
If $\mc= 0$ the spatial parts of the TK and the HM metrics coincide, and this case has already been solved in Section~\ref{ss24XI21.1}.

In the case $\mc\neq0$ and $n\ne 2$ we see from \eqref{HMbeforesplit} that $\hat \gamma(\ntheta)$ has to be a constant. This constant is irrelevant as it can be absorbed into a redefinition of $\beta(r)$, which is why we set it to zero in the following. Equation \eqref{ArrmAttHM} then reduces to
\begin{equation}
    \partial_r ^2 \beta = 0\,.
\end{equation}
Hence, we have that $\beta(r) = r {\tilde c}_1 +{\tilde c}_2 $ which yields
\begin{equation}
\label{HMmnot01}
   V(r,\ntheta, \blue{\theta}^1, \ldots,\blue{\theta}^{n-2}) = r {\tilde c}_1  +{\tilde c}_2 + r 
   c_I \blue{\theta}^I + r \frac{\hat{c}}{2} \blue{|\theta|} ^2 \,.
\end{equation}
Plugging this into the equation \eqref{eqrr}, $A_{rr} = 0$, together with the explicit form of $f(r)$ we obtain
\begin{equation}
    A_{rr}= - \frac{\left(r^n+(-2+n) \mc  \right)}{r^2 (r^n - 2 \mc )}{\tilde c}_2 = 0
    \,,
\end{equation}
from which we conclude that ${\tilde c}_2=0$. The differential equation \eqref{eqthetatheta} $A_{\ntheta \ntheta} = 0$ is then automatically fulfilled. Plugging all of this into
\eqref{thetaIthetaImt} $A_{\blue{\theta}^I \blue{\theta}^I} = 0$ we obtain
\begin{equation}
    A_{\blue{\theta}^I \blue{\theta}^I} = r \hat{c} = 0\,,
\end{equation}
from which follows that $\hat{c} =0$. Hence, we have that
\begin{equation}
\boxed{
\label{HMmnot0final}
   V(r, \ntheta, \blue{\theta}^1, ... \blue{\theta}^{n-2}) = r {\tilde c}_1  + r \sum_{I=1}^{n-2} \left(c_I \blue{\theta}^I \right)\,.
   }
\end{equation}

We finish our analysis with a remark on  space-dimension $n=2$:
      \begin{equation}\label{13I21.1}
        \tilde g = - r^2 dt^2 +  \frac{dr^2 }{r^2   -  2\mc }   + (r^2   -  2\mc  ) d\psi^2
         \,.
      \end{equation}
Setting
$$
 \rho^2 = r^2- 2\mc\,,
$$
one obtains
      \begin{equation}\label{13I21.2}
        \tilde g = - (\rho^2 + 2\mc)  dt^2 +  \frac{d\rho^2 }{\rho^2   + 2\mc }   + \rho^2 d\psi^2
         \,.
      \end{equation}

      When $\mc>0$, regularity at the rotation axis requires that $\psi$ be $2\pi/\sqrt{2\mc}$-periodic.
      To accommodate for this we set
      $$
      \bar \psi =\lambda  \psi
      \,,
      \quad
      \bar t =\lambda  t
      \,,
      \quad
      \bar \rho  =\lambda ^{-1}  \rho
       \,,
      $$
      with $\lambda =  \sqrt{2m_c} $
       chosen so that
       the period of $\bar \psi $ is $2\pi$, leading to
      \begin{equation}\label{13I21.3}
        \tilde g
          = - (\bar \rho^2 + 1 ) d\bar t^2 +  \frac{d\bar \rho^2 }{\bar \rho^2   + 1  }   + \bar \rho^2 d\bar \psi^2
         \,,
      \end{equation}
      which is AdS$_3$ spacetime.

          When $\mc\le 0$, any period of $\psi$ is allowed. In this case, the metric is a BK metric with coordinate mass $-\mc$ if the period of $\psi$ equals $2\pi$.
      As the metric \eqref{13I21.1} is locally BK for $n =2$, the analysis of subsection \ref{s13I22.1} applies.

\section{Transformation behaviour of the mass aspect tensor under conformal transformations}
 \label{A17VIII21.2}
 Consider a metric of the following Fefferman-Graham form
 \begin{equation}
     g = x^{-2} (dx^2 + h_{AB} dx^A dx^B)
     \,,
      \label{30VI21.11}
 \end{equation}
 where the coordinate functions $h_{AB}$ depend upon both $x$ and the local coordinates $x^A$ on the boundary $\{x=0\}$. Here, the index $A$ runs from $1, ... n-1$.
 In the following we will consider $n \geq 3$.
 Let us further write a Taylor expansion
 \begin{equation}
     h_{AB} = (1-k x^2/4)^2\mathring h_{AB}(x^C) + x^n \mu_{AB}(x^C) + o(x^n)
     \,,
 \end{equation}
 where $k$ is a constant. Here, the scalar curvature of $\mathring h_{AB}$ is given by
 \begin{equation}
 \label{riccihd}
     R(\mathring h) = k (n-1)(n-2)\,.
 \end{equation}
 Let $\Rphi =\Rphi (x^C) $ be a function on the boundary, set
 \begin{equation}
  \bar h_{AB}|_{x=0}= \Rphi ^{2} h_{AB}
  \,,
 \end{equation}
 thus \eqref{30VI21.11} takes the form
 \begin{equation}
     g = x^{-2}
     \big(
      dx^2 +  (\Rphi ^{-2} \bar h_{AB} + O(x))
      dx^A dx^B
      \big)
     \,.
      \label{30VI21.12}
 \end{equation}
 We wish to rewrite this as
 \begin{equation}
     g = y^{-2}
     \big(
      dy^2 +  \bar h_{AB}
      d\bar x^A d\bar x^B
      \big)
     \,.
      \label{30VI21.12b}
 \end{equation}
 with
 \begin{equation}
 \label{hbarform}
     \bar h_{AB} = (1- \bar {k}y^2/4)^2 \Rphi ^2 (\bar x^C)\mathring h_{AB} (\bar x^C) + y^n \bar \mu_{AB} (\bar x^C) + o(y^n)
     \,,
 \end{equation}
 where $\bar k \in \{0,\pm 1\}$ {and}
 \begin{equation}
 \label{riccihdb}
    \bar R(\Rphi ^2 \mathring h) = \bar{k} (n-1)(n-2)\,.
 \end{equation}

 By matching powers in Taylor expansions, this is equivalent to finding a function $y$ with a Taylor expansion
 \begin{subequations}
 \label{changecoord}
 \begin{equation}
     y = \Rphi  x \left( 1 + \Rphi _1 x + \Rphi _2 x^2 + \Rphi _3 x^3 + ... \right)
     \label{8VII21.1}
 \end{equation}
 and new boundary coordinates $\bar x^A$ with Taylor expansions
 \begin{equation}
     \bar x^A = x^A+ \Rphi_1^A x + \Rphi_2^A x^2 + \Rphi_3^A x^3+  \Rphi_4^A x^4 ...
     \label{8VII21.3}
 \end{equation}
 \end{subequations}

 In order to determine $y$ as a function of the original coordinates, we note that
 \begin{equation}
     g^{yy}
     \equiv
     g(dy,dy) = y^2
      \qquad
     \Longleftrightarrow
     \qquad
     |d(\red{\log} y)|^2_g= 1
     \,.
 \end{equation}
 This equation says that the integral curves of $d(\red{\log} y)$ are affinely parameterised geodesics. The solutions can be found by shooting geodesics orthogonally, in the metric $x^2g$, to the conformal boundary, with suitable boundary conditions determined by the function $\Rphi $. Hence smooth solutions always exist, which justifies the existence of the expansion \eqref{8VII21.1}. Further, we see that the equation for $y$ can be solved independently of the equation for the $\bar x^A$'s.
 Next, the equation
 \begin{equation}
     0 = g(dy,d\bar x^A)
     \label{8VII21.2}
 \end{equation}
 says that the coordinates $\bar x^A$ are constant along the integral curves of $dy$. So, when a smooth function $y$ is known, one obtains smooth functions $\bar x^A$ by solving \eqref{8VII21.2} along the integral curves of $dy$, which again justifies the expansion
     \eqref{8VII21.3}, and provides a prescription how to find the expansion functions $\psi^A_n$.

 We find that
 \begin{subequations}
 \label{coefflowerorders}
 \begin{equation}
 \Rphi _1 =0 \,,\qquad
  \Rphi _2 = - \frac{\Db_A \Rphi  \Db^A \Rphi }{4 \Rphi^2 }\,, \qquad
 \Rphi _3 =0
 \,,
 \end{equation}
 and
 \begin{align}
     \Rphi^A_1   = 0
  \,, \qquad
  \Rphi^A_2
   =
   -
   \frac 12 \frac{{\Db}^A \Rphi }{\Rphi }\,,
  \qquad
  \Rphi^A_3 = 0\,,
 \end{align}
 \end{subequations}
 where $\Db$ is the covariant derivative of the metric $\mathring h$ and the indices $A$
 are raised with $\mathring h^{A B}$.
 With this it holds that
 \begin{equation}
     g = x^{-2}
     \big(
      dx^2 +  ( h_{AB} + O(x^2))
      dx^A dx^B
      \big)
      =y^{-2}
     \big(
      dy^2 +  ( \bar h_{AB} + O(y)^2)
      d\bar{x}^A d\bar{x}^B
      \big)
 \end{equation}
 However, for general $\mathring h_{AB} $
 it is not possible to bring the $y^2 d\bar{x}^A d\bar{x}^B$ term into the form \eqref{hbarform} by changing coordinates as \eqref{changecoord} after the conformal transformation.
 Indeed, one finds that the terms of order $y^0$ in the coordinate-transformed metric read
 \begin{equation}
 \label{conditionaBC}
     -k \frac{\mathring{h}_{A B}}{2} + \Db_A \Db_B \red{\log}  \Rphi - (\Db_A \red{\log}  \Rphi) (\Db_B \red{\log}  \Rphi)
     + \frac{1}{2} (\Db^C \red{\log}  \Rphi) (\Db_C \red{\log}  \Rphi) \mathring{h}_{A B} \,.
 \end{equation}
 It follows that for \eqref{hbarform} to hold true the following expression must vanish
 \begin{equation}
 \label{condition}
  \frac{\Rphi ^2 }{2} (-\bar{k} \Rphi ^2 +k) \mathring  h_{A B} - \Rphi  \Db_A \Db_B \Rphi
  +2 (\Db_A \Rphi )(\Db_B \Rphi  ) - \frac{1}{2} (\Db^C \Rphi )  (\Db_C \Rphi ) \mathring h_{AB} = 0\,.
 \end{equation}

 The trace of \eqref{condition} is equivalent to the transformation of the Ricci scalar under conformal transformations
 (compare e.g.\ \cite[Appendix D]{Wald:book})
 \begin{equation}
 \label{trafoRicci}
   \bar{R} = \Rphi ^{-2} \left(  R- 2 (n-2) \Db^A \Db_A \red{\log}\Rphi  -  (n-3) (n-2)(\Db^A \red{\log} \Rphi ) (\Db_A \red{\log} \Rphi ) \right)\,.
 \end{equation}
 where we have used the explicit expressions for the Ricci scalars \eqref{riccihd} and \eqref{riccihdb}.

 \subsection{Relation to coordinates in the main text}
  \label{sApp23IX21.1}
 The coordinates $x, x^A$ in \eqref{30VI21.11} are related to the coordinates in the main text $r, x^A$, Section \ref{s3VII21.1} as
 \begin{subequations}
     \label{coordinatetomain}
 \begin{align}
 r =&{} \frac{1}{x} - \frac{k x}{4}
 \,, \\
 x =&{} \frac{2}{r+\sqrt{r^2+k}} = \frac{1}{r}- \frac{k}{4 r^3}+ \frac{k^2}{8 r^5}+O(r^{-7})
 \end{align}
 \end{subequations}
 Under this change of coordinates the line element
 \begin{align}
     g ={}& x^{-2} (dx^2 + \left((1-k x^2/4)^2\mathring h_{AB}(x^C) + x^n \mu_{AB}(x^C) + o(x^n)\right) dx^A dx^B)
      \label{30VI21.11again}
 \end{align}
 changes as
 \begin{subequations}
 \begin{align}
     g ={}& \frac{dr^2}{r^2+k} + r^2 \left(\mathring h_{AB} (x^c) + r^{-n} \mu_{AB}(x^C) + o(r^{-n})\right) dx^A dx^B 
 \end{align}
 \end{subequations}
 Using the relations \eqref{coordinatetomain} we have that
     \begin{align}
         \bar x^A ={}& x^A -
         \frac 12 \frac{\Db^A \Rphi }{\Rphi } x^2  + O(x^4)
         \label{8VII21.1n}
     \end{align}
 which, using
 $$\rphi := \Rphi^{-1}
 \,,
 $$
 translates to
 \begin{align}
     \label{expxAappendix}
    \bar{x}^A = x^A +\frac{\Db^A \psi }{2 \psi r^2}  + O(r^{-4})
         \,,
 \end{align}
 while the expansion in $y$
     \begin{equation}
         y = \Rphi  x \left( 1 -  \frac{{\Db}^A \Rphi {\Db}_A \Rphi}{4 \Rphi^2} x^2 + \mathcal{O}(x^4) \right)
         \label{8VII21.1n2}
     \end{equation}
 leads to an expansion in $\bar{r}$
 \begin{align}
 \bar{r} ={}& \psi r \left(1 +
 \frac{\left({\Db}^A \psi  {\Db}_A \psi + k \psi^2 - \bar{k} \right) }{4 \psi^2 r^2}
 + O(r^{-4}) \right)\,.
  \label{23IX21.1}
 \end{align}
 Using \eqref{trafoRicci} we can also write
 \begin{align}
     \label{exprappendix}
 \bar{r} ={}& \psi r \left(1 +
 \frac{\left(-{\Db}^A {\Db}_A\red{\log}  \psi + (n-2){\Db}^A \red{\log}\psi  {\Db}_A \red{\log}\psi\right) }{2(n-1) r^2}
 + O(r^{-4}) \right)\,.
 \end{align}

 \subsection{\texorpdfstring{$n=3$}{n = 3}}
 For $n=3$, the boundary metric is two-dimensional. Thus we have the relation that
 \begin{equation}
 \label{Riccitensor2d}
     R_{A B} = \frac{R}{2} \mathring {h}_{A B} = k \mathring {h}_{A B} \,, \qquad  \bar{R}_{A B} = \frac{\bar{R}}{2} \Rphi ^2 \mathring {h}_{A B} = \bar{k} \Rphi ^2 \mathring {h}_{A B}\,.
 \end{equation}
 Reconsidering \eqref{condition} and using the expression \eqref{Riccitensor2d} for the Ricci tensor we get
 \begin{equation}
 \label{condition3}
    \bar{R}_{AB} =   R_{AB}- 2 \Db_A \Db_B \red{\log}  \Rphi
    + 2 \Db_A \red{\log}  \Rphi \Db_B \red{\log}  \Rphi- \mathring  h_{A B}\Db^C \red{\log}  \Rphi \Db_C \red{\log}  \Rphi
     \,.
 \end{equation}
 This coincides with the transformation behaviour of the Ricci tensor in two dimensions if and only if
 \begin{equation}
     \mathring h_{AB} \left(\Db^C \Db_C \red{\log}  \Rphi  - \Db^C \red{\log}  \Rphi \Db_C \red{\log}  \Rphi\right)
     =  2 \left( \Db_A \Db_B \red{\log}  \Rphi - \Db_A \red{\log}  \Rphi \Db_B \red{\log}  \Rphi\right)\,.
 \end{equation}
 Using \eqref{condition} the $y^2 d\bar{x}^A d\bar{x}^B$ term can be brought into the correct form.
 For $n=3$ we find that the new mass aspect tensor $\bar{\mu}_{AB}$ takes the form
 \begin{equation}
     \bar{\mu}_{AB} = \frac{\mu_{AB}}{\Rphi } = \rphi \mu_{A B}
     \,,
 \end{equation}
 where $\rphi = \Rphi^{-1}$, which is the variable sometimes used in the main text.
 From
 \begin{equation}
     \mathring h^{AB} \mu_{AB}
     = \Rphi ^3 \bar h_{AB}|_{x=0} \,  \bar \mu_{AB}
     \,,
     \qquad
     \sqrt{\det \mathring h} = \Rphi ^{-2}\sqrt{\det \bar h}|_{x=0}
     \,,\qquad \bar{x}^A = x^A |_{x=0}
     \,,
 \end{equation}
 we conclude that
 \begin{equation}
     \int \mathring h^{AB} \mu_{AB} \sqrt{\det \mathring h} \, dx^1 dx^2
      =
     \int
     \left(
     \Rphi
      \bar h^{AB} \bar \mu_{AB} \sqrt{\det \bar h}
      \right) \big|_{\bar x=0} \, d\bar{x}^1 d\bar{x}^2
      \,,
 \end{equation}
 as well as
 \begin{equation}
     \int
     \left(
      \bar h^{AB} \bar \mu_{AB} \sqrt{\det \bar h}
     \right) \big|_{\bar x=0} \, d\bar{x}^1 d\bar{x}^2
       =
       \int  \Rphi^{-1} \mathring h^{AB} \mu_{AB} \sqrt{\det \mathring h} \, dx^1 dx^2
      \,.
       \label{19IX21.1}
 \end{equation}

 \subsection{\texorpdfstring{$n = 4$}{n = 4}}
  \label{ss20IX21.1}

 We consider now the transformation formulae for the higher-order terms in the expansion of the metric. These are irrelevant for the mass aspect tensor when $n=3$ but become relevant in higher dimensions. Because the powers of $x$ in the expansion of $y$ and $\bar x ^A$ jump by two, and so do the powers in the expansion of the metric up to the mass-aspect level $x^{n-2}$, when $n\ge 4$ the terms of order $x$ in the physical metric remain zero after the change of coordinates. The terms of order $x^2$ change in a non-trivial way, which   is relevant for the mass aspect function in space-dimension $n=4$, and which we determine now.

 The explicit form of the coefficients
 $\Rphi _1, \Rphi _2, \Rphi _3$ and $\Rphi_1^A, \Rphi_2^A, \Rphi_3^A$ in \eqref{coefflowerorders} holds true for arbitrary dimensions $n\geq 3$.
 However,
 to determine the transformation behaviour of the mass aspect in $n=4$ we need to determine further coefficients. We have that
 \begin{equation}
     \Rphi_4^A = - \frac{\bar{k}}{8} \Rphi  \Db^A \Rphi
     + \frac{1}{8 \Rphi ^2}\Db^B \Rphi  \Db^C \Rphi  \Db_C \Db_B x^A \label{psi4A}\,,
     \quad
     \Rphi _4 = - \frac{\bar{k} }{16} \Db_C \Rphi \Db^C \Rphi
     + \frac{1}{16 \Rphi ^4} (\Db_C \Rphi \Db^C \Rphi )^2\,,
 \end{equation}
 where $x^A$ in the first equation  \eqref{psi4A} is treated as a scalar for every $A$.
 In dimension $n= 4$ the coefficients $\Rphi_5^A$ vanish
 \begin{equation}
      \Rphi_5^A =  0\,, \qquad
      \Rphi _5 = 0\,.
 \end{equation}
 This leads to
 \begin{align}
 \label{mufourplusone}
     \Rphi ^2 \bar{\mu}_{AB} ={}& \mu_{AB}+ \frac{1}{16} \mathring h_{AB}\left( {k}^2
     - \bar{k}^2 \Rphi ^4 \right) - \frac{\bar{k}}{4} \Rphi  \Db_A \mathring
     D_B \Rphi
     - \frac{3}{8} \bar{k} \mathring h_{AB} \Db^C \Rphi  \Db_C \Rphi
     + \frac{3}{4} \bar{k} \Db_A \Rphi  \Db_B \Rphi
     \nonumber \\
     &
     - \frac{1}{4 \Rphi ^2} \Db^C \Db_A \Rphi \Db_C \Db_B \Rphi
     + \frac{3}{2 \Rphi ^3} \Db^C \Db_{(A} \Rphi  \Db_{B)} \Rphi  \Db_C \Rphi
     - \frac{1}{4 \Rphi ^3}\mathring h_{AB}
     \Db^C \Rphi  \Db_D  \Db_C \Rphi  \Db^D \Rphi\nonumber \\
     &
     - \frac{1}{2 \Rphi ^3} \Db_A\Db_B \Rphi  \Db^C \Rphi  \Db_C \Rphi
     - \frac{3}{4 \Rphi ^4}\Db_A \Rphi  \Db_B \Rphi \Db^C \Rphi  \Db_C \Rphi
     +\frac{3}{16 \Rphi ^4} (\Db^C \Rphi  \Db_C \Rphi )^2 \mathring h_{AB} \nonumber \\
     &- \frac{1}{4 \Rphi ^2} R_{A C D B} \Db^C \Rphi  \Db^D \Rphi \,.
 \end{align}
 This is best computed using coordinates in which the metric is diagonal, which can always be done locally in dimension three~\cite{DeTurckYang}.
 Here and elsewhere, the symmetrization is defined as $M_{(A B)} = \frac{1}{2} \left(M_{AB}+ M_{BA} \right)$.
 Contracting \eqref{mufourplusone} with $\mathring h^{AB}$ we obtain
 \begin{align}
 \label{trace4p1}
     \Rphi ^2 \bar{\mu}_{AB} \mathring h^{AB} - \mu_A^{~A}={}&  \frac{3}{16} ({k}^2 - \bar{k}^2 \Rphi ^4) - \frac{\bar{k}}{4} \Rphi  \Db^A \Db_A \Rphi
     - \frac{3}{8} \bar{k} (\Db^A \Rphi  \Db_A \Rphi ) - \frac{1}{4 \Rphi ^2} (\Db^A \Db^B \Rphi )(\Db_A \Db_B \Rphi ) \nonumber \\
     & + \frac{3}{4 \Rphi ^3} \Db^A \Rphi  \Db^B \Rphi   \Db_A \Db_B \Rphi - \frac{1}{2\Rphi ^3} (\Db^A \Db_A \Rphi )(\Db^B \Rphi  \Db_B \Rphi ) - \frac{3}{16 \Rphi ^4}(\Db^A \Rphi  \Db_A \Rphi )^2  \nonumber \\
     &+ \frac{1}{4 \Rphi ^2} R_{AB}\Db^A \Rphi  \Db^B \Rphi
     \,.
 \end{align}
 Using the transformation behaviour of the Ricci scalar \eqref{trafoRicci},
 \begin{align}
 \label{trafoRiccin}
   \bar{R}
   &= \Rphi ^{-2} \left(R - \frac{4}{\Rphi }  \Db^A \Db_A \Rphi  + \frac{2}{\Rphi ^2}\Db^A \Rphi  \Db_A \Rphi  \right)
   \,,
 \end{align}
 and the explicit expression of the Ricci scalars \eqref{riccihd} and \eqref{riccihdb},
 \begin{equation}
 \label{Riccinn}
     R=  6 {k}\,, \qquad \bar{R} =  6 \bar{k}\,,
 \end{equation}
 for $n=4$
 we find that
 \begin{equation}
   \bar{k} =  \frac{3 {k} \Rphi ^2 + (\Db^A \Rphi  \Db_A \Rphi ) - 2 \Rphi  \Db^A \Db_A \Rphi }{3 \Rphi ^4}\,,
 \end{equation}
 which yields
 \begin{align}
 \label{trace4p1simp2}
 \Rphi ^2 \bar{\mu}_{AB} \mathring h^{AB} - \mu_A^{~A}
     = {}& - \frac{1}{3 \Rphi ^4} (\Db^A \Rphi  \Db_A \Rphi )^2 - \frac{1}{4 \Rphi ^3}(\Db^A \Rphi  \Db_A \Rphi )(\Db^B \Db_B \Rphi ) + \frac{ (\Db^A \Db_A \Rphi )^2}{12 \Rphi ^2} \nonumber \\
     &  - \frac{1}{4 \Rphi ^2} (\Db^A \Db^B \Rphi )(\Db_A \Db_B \Rphi ) + \frac{3}{4 \Rphi ^3} \Db^A \Rphi  \Db^B \Rphi   \Db_A \Db_B \Rphi   \nonumber \\
     & +\frac{1}{4 \Rphi ^2} (R_{AB} - 2 k \mathring h_{AB})\Db^A \Rphi  \Db^B \Rphi\,.
 \end{align}
 Making the change of variable
 \begin{equation}
     \Rphi  = \exp(u)
 \end{equation}
 (so that $u$ here corresponds to $+\mv/2$ in~\eqref{5VII21.2qw}),
 equation \eqref{trace4p1simp2} becomes
 \begin{align}
     \label{trace4p1simp3}
     \exp(2 u)\bar{\mu}_{AB} \mathring h^{AB} - \mu_A^{~A} = {}& - \frac{1}{12} (\Db^A u \Db_A u)(\Db^B \Db_B u)
     + \frac{1}{12} (\Db^A \Db_A u)(\Db^B \Db_B u)
     \nonumber \\
     &- \frac{1}{4} (\Db^A \Db^B u) (\Db_A \Db_B u) + \frac{1}{4} (\Db^A \Db^B u) (\Db_A u\Db_B u) \nonumber \\
     &+\frac{1}{4} (R_{AB} - 2 k \mathring h_{AB})\Db^A u \Db^B u\,.
 \end{align}
 Using
 \begin{equation}
     \sqrt{\det \mathring h} = \Rphi ^{-3}\sqrt{\det \bar h}|_{x=0}\,, \quad \bar{x}^A = x^A|_{x=0}\,,
 \end{equation}
 and \eqref{trace4p1simp3} we find that
 \begin{equation}
     \int  \left( \bar \mu_{AB} \bar h^{AB} \sqrt{\det \bar h}\right)|_{x=0}\, d\bar x^1 d \bar x^2 d \bar x^3  =
     \int e^{- u}\left( \mu^A_{~A} +C(x^A) \right) \sqrt{\det \mathring h}\, dx^1 dx^2 dx^3
     \,,
 \end{equation}
 with
 \begin{align}
 C(x^A)
 = {}&
      \frac{1}{6} (\Db_A \Db_B \Db^A u) (\Db^B u)
     + \frac{1}{12} \left(R_{AB} - \frac{R}{2}\mathring h_{AB} \right)\Db^A u \Db^B u  \,.
 \end{align}

\bigskip

\noindent{\sc Acknowledgements:}
We are grateful to Michael Anderson, Greg Galloway, Colin Guillarmou, Jack Lee, Yanyan Li,  Rafe Mazzeo, Bruno de Oliveira and Laurent V\'eron for useful discussions.
PTC was supported in part by
the Austrian Science Fund (FWF) under project  P29517-N27  and by
the Polish National Center of Science (NCN) 2016/21/B/ST1/00940. Part of this research was performed while PTC was visiting the Institute for Pure and Applied Mathematics (IPAM) at UCLA, which is supported by the National Science Foundation (Grant No. DMS- DMS-1925919).
ED was supported by the grant ANR-17-CE40-0034   (project CCEM) and  by the ANR EINSTEIN-PPF of the French National Research Agency.
RW was supported by the Austrian Science Fund FWF under the Doctoral Program W1252-N27 Particles and Interactions.
RW thanks the University of Amsterdam for hospitality during the completion of this
project.

%

\bibliographystyle{amsplain}
\bibliography{%
../../references/reffile,%
../../references/newbiblio,%
../../references/hip_bib,%
../../references/newbiblio2,%
../../references/bibl,%
../../references/howard,%
../../references/bartnik,%
../../references/myGR,%
../../references/newbib,%
../../references/Energy,%
../../references/netbiblio,%
../../references/prop,%
../../references/dp-BAMS,%
../../references/PDE,%
CDW-minimal%
}

\end{document}